\documentclass{article}
\usepackage{epsfig,epic,eepic,latexsym, amssymb}
\input xy

\xyoption{all}
\xyoption{all}
\makeindex

\def\pf{PROOF:}
\def    \QED    {\hfill\hbox{\hskip 4pt
                \vrule width 5pt height 6pt depth 1.5pt}}
\def\epf{\QED\\}

\newcommand{\rit}{\mbox{$\hbox{\it I\hskip -2pt R}$}}
\newcommand{\nit}{\mbox{$\hbox{\it I\hskip -2pt N}$}}

\newcommand{\V}{ {\cal V} } 
\newcommand{\R}{\bar{\rit}_+}
\newcommand{\B}{{\bf Bool }}
\newcommand{\p}{ {\cal P} }
\newcommand{\Q}{ {\cal Q} }
\newcommand{\Pz}{  {\cal P}_0 }
\newcommand{\Pu}{ {\cal P}_1 }
\newcommand{\Pd}{ {\cal P}_2 }
\newcommand{\Fp}{Flat_{\cal P}}
\newcommand{\FPu}{ Flat_{{\cal P}_1}  }
\newcommand{\FPd}{ Flat_{{\cal P}_2}  }
\newcommand{\FPz}{ Flat_{{\cal P}_0}  }

\newcommand{\FFu} { WFil }
\newcommand{\FFd} { FFil }
\newcommand{\CFFu}{ CWFil}
\newcommand{\CFFd}{ CFFil}
\newcommand{\F}{{\cal F}}
\newcommand{\Fs}{{\cal F}^s}
\newcommand{\Ba}{\Gamma}
\newcommand{\Bas}{\Gamma^s}
\newcommand{\I}{{\cal I}}
\newcommand{\Mm}{ M^- }
\newcommand{\Mp}{ M^+ }
\newcommand{\Ml}{ M^l }
\newcommand{\Mr}{ M^r }
\newcommand{\raM}{ \tilde{M} }
\newcommand{\Liminf}{ lim^- }
\newcommand{\Limsup}{ lim^+ }
\newcommand{\phicoc}{\phi\mbox{-}\mathit{Cocts}}

\newtheorem{theorem}{Theorem}[section]
\newtheorem{definition}[theorem]{Definition}
\newtheorem{proposition}[theorem]{Proposition}
\newtheorem{lemma}[theorem]{Lemma}

\newtheorem{corollary}[theorem]{Corollary}
\newtheorem{remark}[theorem]{Remark}

\newtheorem{fact}[theorem]{} 
\newtheorem{point}[theorem]{}

\title{Flatness, preorders and general metric spaces}
\author{Vincent Schmitt}

\begin{document}
\maketitle                      

\begin{abstract}
This paper studies a general notion 
of flatness in the enriched 
context: $\p$-flatness where the parameter 
$\p$ stands for a class of presheaves.
One obtains a completion of a category $A$ by
considering the category $\Fp(A)$ of $\p$-flat 
presheaves over $A$.
This completion is related to the free cocompletion
under a class of colimits defined by Kelly.
For a category $A$, for $\p = \Pz$ the class of all 
presheaves, $\FPz(A)$ is the Cauchy-completion 
of $A$. Two classes  $\Pu$ and $\Pd$
of interest for general metric spaces are considered.
The $\Pu$- and $\Pd$- 
flatness are investigated and the associated completions are
characterized for general metric spaces (enrichments over $\R$) 
and preorders (enrichments over $\B$).
We get this way two non-symmetric completions for metric spaces
and retrieve the ideal completion for preorders. 
\end{abstract}

\begin{section}{Introduction}
In \cite{Law73} Lawvere showed amongst other results
that enriched category theory  
was a suitable unifying framework for metric spaces and 
partial orders.
He proved in particular the following.
Preorders and their morphisms as well as general metric
spaces with non-increasing maps occur as categories
and functors enriched over closed monoidal categories.
The base category is $\B$ for preorders
and $\R$ for general metric spaces.
A categorical completion of enrichments may be defined 
so that for the base category $\V = \R$ it amounts to the 
completion \`a la Cauchy of metric spaces 
whereas for $\V = \B$ it
corresponds to the Dedekind-Mac Neille completion of preorders.
Lawvere's categorical completion was therefore just named
Cauchy-completion.\\

Following the spirit of Lawvere's work, one may wonder 
what more theory common to metric spaces and 
preorders may be developed at the categorical level?
The present paper tackles the following problem.
It is known that partial orders admits various 
completions:
\begin{itemize}
\item the Dedekind-Mac Neille completion,
\item the downward completion,
\item the ``algebraic'' or ``ideal'' completion,
\item ...
\end{itemize}
The terminology may vary for the last completions,
but it is clear what they are once said that
\begin{itemize}
\item the Dedekind-Mac Neille completion is defined 
in terms of maximal cuts;
\item the downward completion is in terms of downward closed 
subsets;
\item the algebraic one is in terms of non-empty directed 
down-sets (sometimes called ``ideals'' but we shall avoid 
this confusing terminology).
\end{itemize}
Quite natural questions are whether all these
completions may be described in terms of enrichments and 
if so what they correspond to for metric spaces.\\ 

The answer to this requires a general notion of 
flatness. Flatness in the enriched 
context is already treated in \cite{Kel82-2}, \cite{BQR98}. 
In the last paper the definition ``filtered weights'' relies on 
left Kan extensions preserving certain limits, which is similar 
to the flatness defined in this paper. Nevertheless both these 
works focus on the case when the base $\V$ is locally presentable 
and we shall avoid here such a restriction on $\V$.
Also Street showed 
in \cite{Str83} that the weights of absolute colimits  
form an important class of presheaves related to the 
Cauchy-completion. 
We shall see that this class may be defined in terms of
flatness.
We propose the following definition.
Given a class $\p$ of indexes, 
a presheaf $F:A^{op} \rightarrow \V$ is called ``$\p$-flat'' 
if its left Kan extension along $Y$ preserves all the limits in 
$[A,\V]$ with indexes in $\p$.\\

We have established the following results.
Let us write $\Fp$ for the class of $\p$-flat presheaves.
For any category $A$, 
the full subcategory $\Fp(A)$ of $[A^{op},\V]$ 
with objects $\p$-flat presheaves is its free 
$\Fp$-cocompletion in the sense of \cite{Kel82}.
Calling simply $\Fp(A)$ the ``$\p$-completion'' of $A$,
one obtains therefore a family of completions for 
categories with parameter a family $\p$ of presheaves.
First the free cocompletion of categories is just the
$\p$-completion for $\p$ the empty class of presheaves. 
On the other hand the Cauchy-completion is shown to be the
$\Pz$-completion for $\Pz$ the whole class
of presheaves. Focusing on metric spaces 
we found two more notions of flatness of particular interest.
We define the families
\begin{itemize}
\item $\Pu$ of presheaves over 
the empty category or the unit category $I$ (with one point $*$,
and $I(*,*) = I$) 
\item $\Pd$ of presheaves on categories with finite 
number of objects.
\end{itemize}
In the context $\V=\R$, one may express the $\Pu$-
and $\Pd$- completions in terms of filters
on the metric spaces. This generalizes the fact
that minimal Cauchy filters on a general metric 
space are in one-to-one
correspondence with left adjoint modules on the associated 
category.
We call the filters corresponding to the $\Pu$-flat and
$\Pd$-flat presheaves respectively
{\em weakly flat} and {\em flat}. 
To sum up, let us say that:
\begin{itemize}
\item with the right notion of morphisms,
weakly flat filters, respectively flat filters, occur as  
``non-empty colimits'', respectively ``non-empty filtered
colimits'' of the so-called forward Cauchy sequences.
These sequences were introduced in the literature
as a generalization of Cauchy sequences in non-symmetric
spaces \cite{Sun95}.
\item Cauchy filters are flat, and when the pseudo metric 
is symmetric, flat filters are Cauchy.
\end{itemize}
Eventually, we show that one can
forget category theory and
describe the $\Pu$- and $\Pd$- completions of non-symmetric 
metric spaces in pure topological/metric terms.
The $\Pd$-completion of a symmetric space
amounts to its Cauchy-completion but 
the $\Pd$-completion of a non-symmetric
space certainly generally differs from its 
bi-completion \cite{FL82},\cite{Fla92} and \cite{Sch03}.
In the case $\V = \B$, it appears that the $\Pu$-completion
yields a completion defined in terms of non-empty downward 
subsets whereas the $\Pd$-completion
is the algebraic completion.
For the applications we tried to use as much as possible 
categorical techniques. 
To this respect the only result that seems
not related to category theory
is the characterization of weakly flat/flat filters
in terms of forward Cauchy sequences.\\

This work relies much on the indexed limits/colimits 
computation \`a la Kelly. We adopt the notation and
pick up many results from \cite{Kel82}. We also use
also a little of the 2-categorical theory of enriched modules, 
our references for it are \cite{StWa78}, \cite{BCSW83}, and 
more recently \cite{DaSt97}.
The author has been also much inspired by \cite{BvBR98} and
\cite{Vic}. Both of these works study metric spaces and 
partial orders as enrichments and define completions
by considering ordinary colimits in the presheaf 
categories.\\

The paper is organized as follows.
Section \ref{flatness} treats the notion
of flatness in the enriched context.
Sections \ref{gms} and \ref{preos} are devoted
to applications, respectively to
general metric spaces($\V = \R$) and 
to general preorders ($\V=\B$).
\end{section}

\begin{section}{Flatness}
\label{flatness}
This section treats briefly flatness in the enriched 
context. A generic notion of $\p$-flat presheaf  
where $\p$ stands for a class of indexes
is investigated.
A completion in terms of $\p$-flat presheaves, the $\p$-completion, 
is defined. It is shown to coincide with 
the free $\Fp$-cocompletion in the sense of \cite{Kel82}
where $\Fp$ denotes the class of $\p$-flat presheaves.
If $\p = \emptyset$ the $\emptyset$-completion
is just the free-cocompletion of categories whereas 
for $\Pz$ the class of all presheaves
the $\Pz$-completion amounts to the Cauchy-completion.
Actually $\Pz$-flat presheaves are exactly the presheaves 
which are left adjoint modules. 
We define then two more notions of flatness 
associated with classes $\Pu$ and $\Pd$ of 
presheaves. Their relevance will appear
with the applications in the next sections.\\

We shall consider in this section a symmetric monoidal
complete closed $\V$. For a matter of consistency 
all the $\V$-categories considered are 
by default small, i.e. they have small sets of objects.
We shall precise ``large'' when a category may not be small.
Large $\V$-categories that we shall use
are the category $\V$ itself,
the presheaf categories $[A^{op},\V]$ for small $A$'s
and the categories $PA$ of accessible presheaves 
$A^{op} \rightarrow \V$
for large $\V$-categories $A$.
Indexes of limits and colimits will be also 
considered small.\\

Given a class of presheaves $\phi$, a $\phi$-limit 
(respectively. $\phi$-colimit) is a limit (respectively. colimit) with
index in $\phi$. A functor is $\phi$-continuous  
(respectively. $\phi$-cocontinuous) if and only if 
it preserves $\phi$-limits (respectively. $\phi$-colimits).

\begin{definition}[\p-flatness]
Given a family $\p$ of indexes, 
a presheaf $F: A^{op} \rightarrow \V$ is said
{\em $\p$-flat} when its left Kan extension along $Y$,
$-*F : [A,\V] \rightarrow \V$ preserves all $\p$-limits. 
$\Fp$ will denote the family of all $\p$-flat presheaves,
and for any $\V$-category $A$, $\Fp(A)$ will denote the full 
subcategory of $[A^{op},\V]$ with objects 
$\p$-flat presheaves.
\end{definition}
For any family $\p$ of indexes, 
representables are $\p$-flat since for
any $A(-,a)$, $Lan_Y(A(-,a))$ is the evaluation in $a$ 
that is cocontinuous.\\ 

Since limits and colimits in functor categories
are pointwise, we remind that given functors 
$F: A^{op} \rightarrow \V$, 
$P: K \rightarrow \V$ and 
$G: A \otimes K \rightarrow \V$
equivalent to $G': A \rightarrow [K, \V]$
and also to   $G'': K \rightarrow [A, \V]$, 
one has $(F*G')k \cong F*(G''k)$ and $\{P, G''\}a \cong \{ P, G'a \}$.
Further on we shall use quite freely these isomorphisms.\\

\begin{lemma}
For any $\p$-flat $G: A^{op} \rightarrow \V$ and any functor
$H: A \rightarrow [C^{op},\V]$ with values $\p$-flat functors,
the colimit $G*H$ is again  $\p$-flat.       
\end{lemma}
\pf
Consider $F:K \rightarrow \V$ in $\p$ and 
$L: K \rightarrow [C, \V]$. One has the successive
isomorphisms:
\begin{tabbing}
$Lan_Y(G*H)(\{F,L\})$ \=$\cong$ \=$\{F,L\} * (G* H)$\\
\>$\cong$ 
\>$(G * H) * \{ F,L \}$\\
\>$\cong$ 
\>$G * (H-* \{ F,L \})$ (
\cite{Kel82}, (3.23) ``continuity of a colimit in its index'')\\
\>$\cong$ 
\>$G * (\{ F, H- * L- \})$ 
\=(since for all $a$, $Ha*-$ preserves $\p$-limits)\\
\>$\cong$ 
\>$\{ F, G * (H- * L-) \}$
(since $G*-$ preserves  $\p$-limits)\\
\>$\cong$ 
\>$\{ F, (G * H)*L- \})$ (\cite{Kel82}, (3.23))\\
\>$\cong$ 
\>$\{ F, L- *(G * H) \})$\\
\>$\cong$ 
\>$\{ F, Lan_Y(G * H) \circ L  \})$\\
\end{tabbing}
The resulting isomorphism 
$Lan_Y(G*H)(\{F,L\}) \cong \{ F, Lan_Y(G * H) \circ L  \})$
corresponds actually to the preservation of $\{F, L \}$ by $Lan_Y(G*H)$.\\

{\it Note to the referee: this part of the proof may be omitted}\\
To check this last point, one may consider the following natural isomorphisms:
$$[C,\V](\gamma, \{F,L\}) 
\cong^{\phi^1_\gamma} [K, \V](F,[C, \V](\gamma,L-));$$
$$[A,\V](\tau, H?*\{F,L\}) 
\cong^{\phi^2_{\tau}} [K, \V](F,[A, \V](\tau,H?*L-))$$ 
that 
exhibits $H-*\{F,L\}$ as the limit $\{F, H-*L-\}$;\\
$$[v, G*(H?*\{F,L\})] 
\cong^{\phi^3_v} [K, \V](F,[v,G*(H?*L-)])$$
that exhibits $G*(H-*\{F,L\})$ $\cong$ $G* \{ F , H-*L- \}$ 
as the limit $\{F, G*(H-*L-)\}$.\\

Now the commutation of square $(I)$ on the diagram below
corresponds to the preservation of $\{F,L\}$ by
$H?*-: [C, \V] \rightarrow [A,\V]$.
Also the commutation of square $(II)$ is
the preservation of $\{F,H-*L-\} \cong H-*\{F,L\}$
by $G*-: [A, \V] \rightarrow \V$.
So eventually the outer square commutes which is
the preservation of $\{F,L\}$ by $G*(H-*-)$.

$\xymatrix{
[C,\V](\gamma, \{F,L\}) 
\ar@{}[r]^{ \cong^{\phi^1_\gamma} } 
\ar[d]_{ H?*- } 
\ar@{}[rd]|{(I)}
&
[K, \V](F,[C, \V](\gamma,L-) )
\ar[d]_{[K,\V](F,H?*-)}
\\
[A,\V](H?*\gamma, H?*\{F,L\}) 
\ar@{}[r]|{ \cong^{\phi^2_{H?*\gamma}} } 
\ar[d]_{G*-} 
\ar@{}[rd]|{(II)}
&
[K, \V](F,[A, \V](H?*\gamma,H?*L-))
\ar[d]_{[K,\V](F,G*-)}
\\
[ G*(H?*\gamma)  , G*(H?*\{F,L\})] 
\ar@{}[r]|{\cong_{\phi^3_{G*(H?*\gamma)  }}} 
& 
[K, \V](F,[ G*(H?*\gamma),G*(H?*L-)])
}$\\

Since the isomorphisms 
$\gamma*(G*H) \cong (G*H)*\gamma \cong G*(H-*\gamma)$ 
are natural in $\gamma$, $- * (G*H)$ also preserves
$\{F, L \}$.
\epf

The above lemma has a few consequences that we shall see now.
First, since by Yoneda for any 
presheaf $F: K \rightarrow \V$,
$F \cong F * Y$, one has  
\begin{proposition}
\label{pflatchar}
For any family $\p$ of indexes,
$F$ is $\p$-flat if and only if it is a
$\Fp$-colimit of representables.
\end{proposition}
This together with the following result \ref{freecoc}
from \cite{Kel82} (see also \cite{AK88})
that relates the closure of categories under $\phi$-colimits
to their free ${\phi}$-cocompletion
will yield a universal completion
for categories in terms of $\p$-flat presheaves
(\ref{fcomp}).\\
 
Remember from \cite{Kel82} that given a family $\phi$ of indexes, 
and a category $A$, the {\em closure of $A$ under $\phi$-colimits}, 
say $\bar{A}$, is defined as the smallest full (replete) subcategory 
of $[A^{op}, \V]$, containing the representables and closed under 
the formation of $\phi$-colimits in $[A^{op},\V]$, which means that
for any $G: K \rightarrow [A^{op},\V]$ taking values 
in $\bar{A}$ and any $F \in \phi$,
$F * G$ is in $\bar{A}$. 
$\bar{A}$ is a full subcategory 
of the $\V$-category $[A^{op},\V]$ and may be not small. 

\begin{theorem}
\label{freecoc} For any family of indexes $\phi$,
for any $A$, the closure $\bar{A}$ of $A$
in $[A^{op},\V]$ under $\phi$-colimits 
constitutes its {\em free ${\phi}$-cocompletion}.   
This means that for any categories $A$:
\begin{itemize}
\item 
$\bar{A}$ is $\phi$-cocomplete;
\item For any possibly large category $B$ 
one has an equivalence 
$Lan_K: [A,B]' \cong \phicoc[\bar{A}, B]$ 
where:
\begin{itemize}
\item $K$ is the full and faithful inclusion 
$A \rightarrow \bar{A}$ (sending any $a \in A$ to $A(-,a)$);
\item $[A,B]'$ stands for the full subcategory of $[A,B]$
of functors admitting a left Kan extension along $K$;
\item $\phicoc[\bar{A}, B]$ is the full subcategory
of $[\bar{A},B]$ 
of $\phi$-cocontinuous functors;  
\item $Lan_K$ stands for the ``left Kan extension functor'', it has 
inverse the restriction to $\phicoc[\bar{A}, B]$ of 
$[K,1]: [\bar{A},B] \rightarrow [A,B]$.
In particular if $B$ is $\phi$-cocomplete then $[A,B]' = [A,B]$.
\end{itemize}
\end{itemize}
\end{theorem}

Thus by \ref{pflatchar},
\begin{theorem}
\label{fcomp}
For any family $\p$ of indexes and 
any category $A$, $\Fp(A)$ is the free 
$\Fp$-cocompletion of $A$.
\end{theorem}
We shall therefore simplify the terminology
and call $\Fp(A)$ the {\em $\p$-completion} of $A$,
for any category $A$, and any family of indexes $\p$.
For $\p = \emptyset$ the empty class of presheaves
all the presheaves are $\emptyset$-flat thus the  
$\emptyset$-completion is just the free cocompletion.
On the other hand let $\Pz$ denote the whole class of 
presheaves, then 
the $\Pz$-completion is the Cauchy-completion.
This is a straightforward consequence of 
\begin{theorem} 
\label{isthattrue}
For a presheaf $F: A^{op} \rightarrow \V$ the following
assertions are equivalent:
\begin{itemize}
\item $(1)$ $F$ is $\Pz$-flat;
\item $(2)$ $F$ as a module $I \rightarrow A$ is a 
left adjoint.
\end{itemize}
\end{theorem}
Before to establish \ref{isthattrue}, we need 
\begin{proposition}
\label{commut}
Let $F: A^{op} \rightarrow \V$,
$P: K \rightarrow \V$ and $G: A \otimes K \rightarrow \V$
equivalent to $G': A \rightarrow [K,\V]$ and also 
to $G'': K \rightarrow [A,\V]$. Then   
$F*-: [A, \V] \rightarrow \V$ preserves the limit 
$\{P , G''\}$ if and only if 
$\{P,-\} : [K, \V] \rightarrow \V$ preserves the
colimit $F* G'$.
\end{proposition}
\pf 
Let $$\xymatrix{ P \ar[r]^-{\eta} & [A, \V](\{ P, G''\},G'' - )   }$$
be the unit of $\{ P, G'' \}$ and 
$$\xymatrix{ F \ar[r]^-{\lambda} & [K, \V](G'-, F*G') }$$
be the unit of $F*G'$.
We need to show that 
$$(1)\;\;\xymatrix{ 
P \ar[r]^-{\eta} & 
[A, \V](\{ P, G''\},G'' -) \ar[r]^-{F*-} &
[F* \{P, G''\}, (F*G''-)] }$$
exhibits $F* \{P, G''\}$ as $\{P, F*G' \}$ if and only if 
$$(2)\;\;\xymatrix{ 
F \ar[r]^-{\lambda} & 
[K, \V](F* G',G' -) \ar[r]^-{\{ P,- \}} & 
[\{P, F*G'\}, \{P, G'- \}]  }$$ exhibits
$\{P, F*G' \}$ as $F* \{P, G''\}$.\\

First note that given $x \in \V$, any natural in $k$, 
$$(1')\;Pk \rightarrow_k [x, (F*G')k]$$
corresponds via Yoneda to a natural in $v$,
$[v,x] \rightarrow_{v} [K,\V](P,[v,(F*G')-])$.
Since $$[K,\V](P,[v,(F*G')-] \cong_v [v,\{P, F*G' \}],$$
it corresponds also to an arrow $$(1'')\;\;x \rightarrow \{P, F*G' \}.$$
Also that $(1')$ exhibits
$x$ as the limit $\{P, F*G' \}$ is equivalent to 
the fact that $(1'')$ is iso.
Analogously any natural in $a$, $$(2')\;\; Fa \rightarrow [\{ P, G''\}a, x]$$ 
corresponds to an arrow $$(2'')\;\;F * \{P, G''  \} \rightarrow x,$$ and
$(2')$ exhibits $x$ as the colimit $F* \{P, G''  \}$ if and only if
$(2'')$ is iso.\\

Now the result follows from the fact the arrow $(1)$ above
corresponds by the bijection $(1')-(1'')$ 
to the same arrow as $(2)$ by $(2')-(2'')$.
To check this last point, use the following sequences of isomorphisms
\begin{tabbing}
\hspace{2cm}\=$[K, \V](P, [F*\{P,G''\},(F*G')- ])$\\
\hspace{1cm}$\cong$\> 
$[K, \V](P, [A^{op}, \V](F,  [\{ P,G'' \}-,(F*G')- ] ))$\\
\hspace{1cm}$\cong$\>
$[K \otimes A^{op}, \V](F\otimes P, [\{ P,G'' \},F*G'])$\\
\hspace{1cm}$\cong$\> 
$[A^{op}, \V](F, [K, \V](P, [ \{P,G''\}-,(F*G')- ] ))$\\
\hspace{1cm}$\cong$\> 
$[A^{op}, \V](F,[\{P,G''\}-, \{P, F*G'\}]$.  
\end{tabbing}
Through this sequence of isomorphisms the natural in $k$ 
$$\xymatrix{ Pk \ar[r]^-{\eta} & [A, \V](\{ P, G''\},G'' k)
\ar[r]^-{F*-} & [F* \{P, G''\}, (F*G')k]}$$ 
corresponds to 
the natural in $a$, $k$
$$\xymatrix{
Fa \otimes Pk \ar[r]^-{ \lambda_{k,a} \otimes \eta_{k,a} } &
[G(k,a), F*G''k] \otimes [\{P,G'a\}, G(k,a)  ] \ar[r]^-{\mu} &
[\{ P, G'a \}, F*G''k ] }$$
where $\mu$ is the composition in $\V$,
the latter one corresponds to the natural in $a$,
$$\xymatrix{ 
Fa \ar[r]^-{\lambda} & [K, \V](G'a, F* G')
\ar[r]^-{\{P,-\}} & [\{P, G'a \}, \{P, F*G'\} ] }.$$
\epf

According to \ref{commut},
\begin{fact}
\label{myremark}
a presheaf $F:A^{op} \rightarrow \V$ 
is $\Pz$-flat if and only if for any $G: A \rightarrow [K,\V]$ 
the colimit $F*G$ is preserved by any representable 
$[K,\V] \rightarrow \V$.
\end{fact}
To prove \ref{isthattrue} we will need to use a bit of 
the 2-categorical machinery developed in \cite{StWa78},
in particular the description of indexed colimits
in terms of right liftings.
It is proved in \cite{Str83}that
\begin{theorem}
\label{Ross}
A module 
$\theta: \xymatrix{A  \ar[r]|{\circ} & B}$
is left adjoint if and only if any colimit indexed by 
$\theta$ is absolute.
\end{theorem}
So \ref{Ross} and \ref{myremark} 
give immediately $(2) \Rightarrow (1)$ 
in \ref{isthattrue}.
Now a minor adaptation of the proof
presented in \cite{Str83} 
will show
$(1) \Rightarrow (2)$ in \ref{isthattrue}.
\begin{proposition}
A presheaf $F: A^{op} \rightarrow \V$ is 
a left adjoint as a module $\xymatrix{I \ar[r]|{\circ} & A}$ if 
the colimit $F \cong F*Y$ is preserved
by any representable $[A^{op},\V] \rightarrow \V$.
\end{proposition}
\pf
That $F \cong F*Y$ amounts to saying that
there is a right lifting of $\V$-modules
as below:
\begin{center}
$(*)$ 
$\xymatrix{
& [A^{op},\V] \ar[ld] \ar@{}[d]|{\Rightarrow} \ar[rd]^{Y^*} &\\ 
I \ar[rr]_{F}&  & A}$
\end{center}
where
\begin{itemize}
\item the unlabeled diagonal is the right adjoint module
$\delta^*$ given by the functor 
$\delta: I \rightarrow [A^{op},\V]$ that sends the one 
point to the presheaf $F$;
\item the horizontal arrow, denoted $F$, is the module
$I \rightarrow A$ corresponding to the presheaf $F$.
\end{itemize}
Recall that any left adjoint module
respects right liftings. So by pasting
$Y_*: A \rightarrow [A^{op}, \V]$ to the 2-cell $(*)$
and since $Y^* \circ Y_* = 1$
one obtains a right lifting
\begin{center}
$(**)$ $\xymatrix{ & A \ar[ld] \ar@{}[d]|{\Rightarrow} \ar[rd]^{1}  &\\ 
I \ar[rr]_{F}&  & A}$
\end{center}
That the left diagonal constitutes a right adjoint
to the module $F$, is equivalent to say that this right
lifting is absolute.
Consider any module $\theta: B \rightarrow [A^{op},\V]$.
It may be decomposed in a product $Y^* h_*$ 
for the functor $Y: [A^{op},\V] \rightarrow P[A^{op},V]$ and 
a functor $h: B \rightarrow [A^{op},\V]$ (sending
any $b$ to $\theta(-,b)$). 
By assumption the colimit $F \cong F*Y$ is preserved 
by any representable presheaf $[A^{op},\V] \rightarrow \V$, thus
because colimits in $P[A^{op},\V]$ are ``pointwise'' (\cite{AK88}),
$Y: [A^{op},\V] \rightarrow P[A^{op},\V]$ also preserves the colimits 
indexed by $F$. That is to say that
$Y^*$ respects the right lifting $(*)$. 
$\theta = h_* Y^*$ also respects $(*)$ (since $h_*$ is left adjoint).
So $(*)$ is absolute as well as $(**)$.
\epf

Let us mention another consequence of \ref{commut}.
\begin{definition}[\Q-coflatness]
Given a family $\Q$ of indexes, 
a presheaf $P: K  \rightarrow \V$ is said 
{\em $\Q$-coflat} when $\{P,-\}: [K,\V] \rightarrow \V$
preserves $\Q$-colimits. Let $Coflat_{\Q}$ denote the family of 
all $\Q$-coflat presheaves. 
\end{definition}
Consider the class $CPSh$ of classes of presheaves.
Then the classes of presheaves are partially ordered by inclusion 
and one has a Galois connection 
$$Flat_{-} \dashv Coflat_{-}: CPSh \rightharpoonup CPSh^{op}.$$  

We shall investigate in this paper two more notions of flatness.
Let us define
\begin{definition}
$\Pu$ is the class of indexes of the form 
$F:K \rightarrow \V$ where $K$ is 
the empty $\V$-category or $K = I$. 
$\Pd$ is the class of indexes $F:K \rightarrow \V$
with $Obj(K)$ finite.
\end{definition}

We shall call {\em conical finite limit} a conical limit
indexed by a finite ordinary category.
From now on we write $A_0$ for the underlying ordinary category
of a $\V$-category $A$.
A minor adaptation of the proof of theorem \cite{Kel82} (3.73) 
as in \cite{Kel82-2} (4.3), shows that
\begin{proposition} 
\label{ect1}
A $\V$-category $A$ is $\Pd$-complete if and only if it 
has all conical finite limits and cotensors. 
Given a $\Pd$-complete $A$, a $\V$-functor $P:A \rightarrow B$ 
is $\Pd$-exact if only if it preserves conical finite limits and cotensors.
\end{proposition}
\pf (Sketch of) It suffices to reuse the argument
developed in the sketch of proof of \cite{Kel82-2} (4.3). 
Remark that if $A$ has conical finite
limits and cotensors, the indexed limit $\{F,G\}$ of any
$F: K \rightarrow \V$ and $G: K \rightarrow A$ with $Obj(K)$ finite
may be computed as the equalizer in $A_0$ of  
$$\xymatrix{ \prod_{k \in K} Fk \pitchfork Gk \ar@<1ex>[r]
\ar@<-1ex>[r] &
\prod_{k,k' \in K} K(k,k') \pitchfork (Fk
\pitchfork Gk')}.$$ Actually all the ordinary limits involved 
in this equalizer, i.e. the two products and the equalizer itself 
are finite and thus conical.  
Also revisiting the sketched proof of theorem (3.73) in \cite{Kel82},
one gets that any functor $H: A \rightarrow B$ preserving
conical finite limits and cotensors, $H$ will preserve the 
above conical equalizer which image in $B$ is then the limit 
$\{F, HG\}$.
\epf

\begin{proposition}
\label{Lprpf}
Given a $\p$-flat presheaf $F: A \rightarrow \V$,
and any functor $G: A \rightarrow B$, the 
left Kan extension of $F$ along $G$ is $\p$-flat.
\end{proposition}
\pf
Given a $\p$-flat $F$, it is a $\Fp$-colimit of representables.
The image by $Lan_G$ of any representable is again 
representable  
(For any $a \in A$, $Lan_G(A(a,-))(b) \cong_b B(G-,b) * A(a,-) 
\cong_b B(Ga,b)$).  
Also the left Kan extension functor 
$Lan_G : [A, \V] \rightarrow [B, \V]$ is cocontinuous as
shown below \ref{LKcoc}, so 
$Lan_G(F)$ is also a $\Fp$-colimit of representables
and thus $\p$-flat according to \ref{pflatchar}.
\epf

\begin{lemma}
\label{LKcoc}
Given any functor $G: A \rightarrow B$,
the left Kan extension functor 
$Lan_G : [A, \V] \rightarrow [B, \V]$ is cocontinuous.
\end{lemma}
Given $J: K^{op} \rightarrow \V$ and $H: K \rightarrow [A,\V]$,
one has the following pointwise computation in $b \in B$, 
\begin{tabbing}
\hspace{1cm}$Lan_G(J*H)(b)$\=$\cong$\=$\tilde{G}b*(J*H)$,\\
\>$\cong$ \>$(J*H) * \tilde{G}b$\\
\>$\cong$ \>$J * (H- * \tilde{G}b)$\\
\>$\cong$ \>$J * (\tilde{G}b * H-)$\\
\>$\cong$ \>$J * ( Lan_G(H-)(b) )$\\
\>$\cong$ \>$( J * Lan_G(H-) ) (b) $.
\end{tabbing}   
Actually the resulting natural isomorphism exhibits 
$Lan_G(J*H)$ as the colimit $J * Lan_G(H-)$.
$Lan_G:[A,\V] \rightarrow [B,\V]$ preserves
 $J* H$ if and only if for any $b \in B$,
$E_b \circ Lan_G: [A,\V] \rightarrow \V$ does. 
But $E_b \circ Lan_G$ 
$\cong$ $Lan_G(-)(b)$ $\cong$ $\tilde{G}(b) * -$ $=$ 
$- * \tilde{G}(b)$ that is known cocontinuous. 
\epf


The rest of the paper treats
notions of flatness for enrichments over 
particular bases namely $\V= \B$ and $\V = \R$.
An important point to make is that 
for both these cases the base $\V$ is 
small and thus is necessarily a preorder 
(see \cite{Bor94} prop. 2.7.1 p.59). 
In the case of a small $\V$, for any small
$\V$-category $A$, the presheaf category $[A,\V]$ remains
small and so does $\Fp(A)$
for any family $\p$ of presheaves.
Still in this case, if $A$ is $\p$-complete then
it is a retract of $\Fp(A)$ (i.e. the inclusion 
$A \hookrightarrow \Fp(A)$ is a split monic) 
but it is generally NOT
isomorphic to $A$.
\end{section}

\begin{section}{The case $\V = \R$.}
\label{gms}
This section treats flatness in the context 
of general metric spaces.
First we come back quickly in \ref{Laws} on Lawvere's 
Cauchy-completion of general metric spaces.
In \ref{modfil}, the existing correspondence 
between Cauchy filters and left adjoint modules
is extended:
the ordinary category of $\FPu$-modules is reflective 
in a category of particular filters,
the so called {\em weakly-flat} ones, with 
reverse inclusion ordering.
By considering the category of fractions
induced by this full reflection one defines a 
notion of morphisms of weakly flat filters 
that yields an enriched equivalence with the 
categories of $\FPu$-modules.
This equivalence restricts to the category
of $\FPd$-modules on one side and on the other side
to a full subcategory of filters, the so-called {\em flat} 
ones. Also in the symmetric case flat filters
are Cauchy and one retrieves via the latter
equivalence the well known one-to-one correspondence 
between Cauchy filters and left adjoint modules.
Weakly flat and flat filters are then related
to forward Cauchy sequences.
These sequences were introduced in the literature
as a generalization of Cauchy sequences \cite{Sun95}. 
They are relevant as with the right notion of morphisms, 
both weakly flat and flat filters occur as canonical 
colimits of functors with values these forward Cauchy 
sequences. In \ref{mcomp}, the $\Pu$-
and $\Pd$-cocompletions of general metric spaces 
are defined and ``internally'' described
in pure metric/topological terms by means
of the previous filters.
A few examples of these completions follow.

\subsection{Lawvere's completion}
\label{Laws}
Let us recall a few results that are from \cite{Law73} or 
belong to folklore.\\

$\R$ stands for the monoidal closed category with:
\begin{itemize}
\item objects: positive reals and $+\infty$;
\item arrows: the reverse ordering, $x \rightarrow y$ if and 
only if $x \geq y$;
\item tensor: the addition (with $+\infty + x = x + +\infty =  +\infty$); 
\item unit: 0.
\end{itemize}
For any pair $x,y$ of objects in $\R$, the exponential object 
$[x,y]$ is $max\{ y - x, 0\}$.\\

A {\em $\R$-category} $A$ corresponds to a {\em general metric space}.
It consists in a set of {\em objects} or {\em elements}, $Obj(A)$ 
(sometimes just denoted $A$) together with a 
map $A(-,-):  Obj(A) \times  Obj(A) \rightarrow \R$
that satisfies:
\begin{itemize}
\item for all $x,y,z \in Obj(A)$, $A(y,z) + A(x,y) \geq A(x,z)$; 
\item for all $x \in Obj(A)$, $0 \geq A(x,x)$.
\end{itemize}
A {\em $\R$-functor} $F: A \rightarrow B$
corresponds to a non-expansive map $F: Obj(A) \rightarrow Obj(B)$, i.e.
for all $x,y \in Obj(A)$, $A(x,y) \geq B(F(x), F(y))$. 
A {\em $\R$-natural transformation} $F \Rightarrow G: A \rightarrow B$
corresponds to the fact that
for all $x \in  Obj(A)$, $0 \geq B(F(x), G(x))$.
A $\R$-module $M: \xymatrix{I \ar[r]|{\circ} & A}$ - or left module
on $A$ - is a map 
$Obj(A) \rightarrow \R$ such that for all $x,y \in A$,
$M(y) + A(x,y) \geq M(x)$.
Dually a $\R$-module $N: \xymatrix{A \ar[r]|{\circ} &  I}$ - or right module 
on $A$ - is a map 
$Obj(A) \rightarrow \R$ such that for all $x,y \in A$,
$A(x,y) + N(x) \geq N(y)$.
The presheaf category $[A^{op},\R]$ has homsets
given by $[A^{op},\R](M,N) = \bigvee_{x \in A}[M(x),N(x)]$.
Its underlying category is a partial order
with arrows given by the pointwise reverse ordering
$M \Rightarrow N$ if and only if $\forall x \in A$,
$M(x) \geq N(x)$.
The composition of left and right modules is as follows.
Given $\xymatrix{ I \ar[r]|{\circ}^M & A \ar[r]|{\circ}^N & I}$, the
composite $N * M$ is $\bigwedge_{x \in A} M(x) + N(x)$.
For such $M$ and $N$,
$M$ is left adjoint to $N$
if and only if:
\begin{itemize}
\item $(1)$ $0 \geq N * M$; 
\item $(2)$ for all $x,y \in A$, $N(y) + M(x)  \geq A(x,y)$.
\end{itemize}

The key point for the Cauchy-completion of general metric 
spaces is that for a general metric space $A$
there is a one-to-one correspondence between left adjoint
modules on $A$ and minimal Cauchy filters on $A$.
From this observation mainly, one gets that the full 
subcategory of $[A^{op},\R]$
with objects left adjoint modules is isomorphic
the completion ``\`a la Cauchy'' of $A$, that is 
its Cauchy-completion if $A$ is a metric space
or more generally its bi-completion if the space
is not symmetric (see \cite{FL82} and \cite{Fla92}
or \cite{Sch03} for the connection with Lawvere's work).  
As this is the starting point of our investigation,
we recall briefly this correspondence.\\ 

Let $A$ stand for a general metric space.
\begin{definition}
A filter $\F$ on $A$ is {\em Cauchy} if 
and only if for any $\epsilon > 0$, there exists an $f \in \F$
such that for any elements
$x, y$ of $f$, $A(x,y) \leq \epsilon$ or equivalently when: 
$$\bigwedge_{f \in \F} \bigvee_{x,y \in f} A(x, y)= 0.$$
\end{definition}

\begin{definition}
For any left adjoint module $M$ on $A$, with right adjoint $\raM$
one defines $\Bas(M)$ as the subset of 
$\wp(A)$:
$\{ \Bas(M)(\epsilon) \mid \epsilon \in ]0, + \infty] \}$,
where $\Bas(M)(\epsilon)$ denotes the set $\{ x \in A \mid 
M(x) + \raM(x) \leq \epsilon\}$.  
\end{definition}

For any left adjoint module $M$ on $A$, $\Bas(M)$ is a Cauchy 
basis. The filter that it generates, that we denote $\Fs(M)$, is a 
minimal Cauchy filter. The map $M \mapsto \Fs(M)$ defines a
bijection between left adjoint 
modules $\xymatrix{ I \ar[r]|{\circ} & A}$
and minimal Cauchy filters on $A$.
Actually one may check the following points (see for example \cite{Sch03}).
To any Cauchy filter $\F$ one may associate 
a left adjoint module $\Ml(\F)$ defined by
$$x \mapsto \bigwedge_{f \in \F} \bigvee_{y \in f} A(x,y) 
         = \bigvee_{f \in \F} \bigwedge_{y \in f} A(x,y).$$
$\Ml(\F)$ has right adjoint $\Mr(\F)$ given by the map
$$x \mapsto \bigwedge_{f \in \F} \bigvee_{y \in f} A(y,x)= 
\bigvee_{f \in \F} \bigwedge_{y \in f} A(y,x).$$
For any left adjoint module $M$ on $A$, 
$\Ml(\Fs(M)) = M$ and
for any Cauchy filter $\F$ on $A$, 
$\Fs(\Ml(\F))$ is the only
minimal Cauchy filter contained in $\F$.\\

{\it Note to the referee - to be omitted}\\
For what it is worth. It is well known that any Cauchy filter 
contains only one minimal Cauchy filter. But I don't know from 
the literature - apart from \cite{Sch03} - any explicit proof
that for any Cauchy filter $\F$ on $A$, 
$\Fs(\Ml(\F))$ is the only minimal Cauchy filter contained 
in $\F$. So here are two key points to retrieve quickly that result
once you suppose that for all left adjoint module
$M$, $M = \Ml \circ \Fs (M)$.\\
$(1)$ For any $\F$ Cauchy, one may check that $\F \supseteq \Fs(\Ml(\F))$ 
using the definition 
$\Ml(\F)(x) = \bigwedge_{f \in \F} \bigvee_{y \in f} A(x,y)$.\\
$(2)$ Also for Cauchy filters $\F_1$ and $\F_2$, if
 $\F_1 \supseteq \F_2$ then 
$\Ml(\F_2) \Rightarrow \Ml(\F_1)$ by using the definition
$\Ml(\F)(x) = \bigwedge_{f \in \F} \bigvee_{y \in f} A(x,y)$
and $\Ml(\F_1) \Rightarrow \Ml(\F_2)$ by using 
the definition $\bigvee_{f \in \F} \bigwedge_{y \in f} A(x,y)$,
so eventually  $\Ml(\F_1) = \Ml(\F_2)$.

\subsection{Modules and Filters}
\label{modfil}
A  natural question is whether the previous 
correspondence left adjoint modules / Cauchy filters
may be extended to a class of $\p$-flat modules. We shall 
show that this is the case for $\p$ = $\Pu$ and $\Pd$.\\

Let $A$ denote from now on a general metric space.\\

Let us give an explicit definition of those $\Pu$-flat and
$\Pd$-flat modules.
We shall recall first a few technical points.
For the assertions \ref{cotpoint}, \ref{precon} and 
\ref{colimpre} below, $\V$ denotes a complete monoidal closed $\V$.
Remember that cotensors are defined pointwise in functor
categories. In particular 
\begin{fact}
\label{cotpoint}
Any presheaf $\V$-category $[A,\V]$ is cotensored: for any 
presheaf $N$,
$v \pitchfork N$ is the composite 
$\xymatrix{ A \ar[r]^{N} & \V \ar[r]^{[v,-]} & \V }$.
\end{fact}
Also for functor between cocomplete categories
the preservation of conical colimits amounts
to the preservation of ordinary colimits. Precisely
one may check:
\begin{fact}
\label{precon}
Given a $\V$-functor $T:A \rightarrow B$ 
with underlying ordinary functor
$T_0: A_0 \rightarrow B_0$ 
and an ordinary functor $P: J \rightarrow A_0$,
if the conical limits of $P$ and of $T_0P$ exist
and $T_0$ preserves the ordinary limit of $P$,
then $T$ preserves the conical limit of $P$.
\end{fact}
Eventually the preservation of limits/colimits is simple
in the case $\V = \R$ since
\begin{fact}
\label{colimpre}
If the base category $\V$ is a preorder, then
given a presheaf $F: A^{op} \rightarrow \V$ and
a functor $G: A \rightarrow B$ such that $F*G$ exists
and $H: B \rightarrow C$ then
$H$ preserves $F*G$ if and
only if $F*(GH)$ exists and $H(F*G) \cong F*(GH)$.
\end{fact}

According to the three previous point  
one gets
\begin{fact}
Let $M: \xymatrix{I \ar[r]|{\circ} & A}$ be a left module.
\begin{itemize}
\item $-*M: [A,\R] \rightarrow \R$ preserves 
the unique conical limit with index with empty
domain if and only if the underlying ordinary functor 
preserves the terminal object
i.e. $0 * M = 0$ if and only if 
$$(1)\;\;\bigwedge_{x \in A} M(x) = 0.$$  
\item $-*M$ preserves conical finite limits if and only\\ 
$(2)$ For any finite family of right modules $N_i: \xymatrix{A \ar[r]|{\circ} & I}$, $i \in I$,
$$\bigwedge_{x \in A}( M(x) + \bigvee_{i \in I} N_i(x) )
= \bigvee_{i \in I}( \bigwedge_{x \in A} M(x) + N_i(x)  );$$
\item $-*M$ preserves cotensors if and only if\\
$(3)$ For any $v \in \R$ and 
any right module $N: \xymatrix{A \ar[r]|{\circ} & I}$ , 
$$\bigwedge_{x \in A}( M(x) + [v,N(x)] ) = [v, \bigwedge_{x \in A}( M(x) + N(x) )].$$
\end{itemize}
\end{fact}
So $\Pu$-flat modules are those satisfying $(1)$ and $(3)$ above,
and $\Pd$-flat modules are those satisfying $(2)$ and $(3)$.\\  

It is convenient to introduce now the following notations.
\begin{definition}
Given a filter $\F$ on $A$ and a map $f: Obj(A) \rightarrow Obj(\R)$,
$\Limsup_{x \in \F} f(x)$ or simply $\Limsup_{\F} f$ denotes 
$\bigwedge_{f \in \F} \bigvee_{x \in f} f(x)$. Also
 $\Liminf_{x \in \F} f(x)$ or $\Liminf_{\F} f$ will stand for
$\bigvee_{f \in \F} \bigwedge_{x \in f} f(x)$.
\end{definition}

From the correspondence Cauchy filters/left adjoint modules,
we know two operators that associate filters to modules.
\begin{definition}
Given any filter $\F$ on $A$, we define the following $\R$-valued maps
on objects of $A$: 
\begin{tabbing}
\hspace{1cm}\=$\Mm(\F): x \mapsto \Liminf_{\F} A(x,-) 
= \bigvee_{f \in \F} \bigwedge_{y \in f} A(x,y)$,\\
\>$\Mp(\F): x \mapsto  \Limsup_{\F} A(x,-) = \bigwedge_{f \in \F} \bigvee_{y \in f} A(x,y)$.\\
\end{tabbing}
\end{definition}
For any filter $\F$ on $A$, one has $\Mm(\F) \leq \Mp(\F)$, and if $\F$
is Cauchy then $\Mm(\F) = \Mp(\F)$. 

\begin{fact}
Given any filter $\F$ on $A$, the map $\Mp(\F)$ defines a module 
$\xymatrix{I \ar[r]|{\circ} & A}$.
\end{fact}
\pf
One has to show that for all $x,y \in A$, $\Mp(\F)(x) + A(y,x) \geq \Mp(\F) (y)$.
For all $x,y \in A$, 
\begin{tabbing}
$\Mp(\F)(x) + A(y,x)$ \=$=$ \=$( \bigwedge_{f \in \F} \bigvee_{z \in f} A(x,z) ) + A(y,x)$\\
\>$=$ \>$\bigwedge_{f \in \F} ( ( \bigvee_{z \in f}  A(x,z))  + A(y,x)   )$\\
\>$\geq$\>$\bigwedge_{f \in \F} \bigvee_{z \in f} ( A(x,z) + A(y,x) )$\\
\>$\geq$\>$\bigwedge_{f \in \F} \bigvee_{z \in f} A(y,z)$\\
\>$=$ $\Mp(\F)(y)$.
\end{tabbing}
\epf

\begin{fact}
Given any filter $\F$ on $A$,
the map $x \mapsto \Mm(\F)(x)$ defines a module 
$\xymatrix{I \ar[r]|{\circ} & A}$.
\end{fact}
\pf
One has to show that for all $x,y \in A$, 
$\Mm(\F)(x) + A(y,x) \geq \Mm(\F)(y)$.
For all $x,y \in A$, 
\begin{tabbing}
$\Mm(\F)(x) + A(y,x)$ \=$=$ \=$( \bigvee_{f \in \F} \bigwedge_{z \in f}  A(x,z) ) + A(y,x)$\\
\>$\geq$ \>$\bigvee_{f \in \F} ( (\bigwedge_{z \in f} A(x,z)) + A(y,x)   )$\\
\>$=$\>$\bigvee_{f \in \F} \bigwedge_{z \in f} (A(x,z) + A(y,x))$\\
\>$\geq$\>$\bigvee_{f \in \F} \bigwedge_{z \in f} A(y,z)$\\
\>$=$ $\Mm(\F)(y)$.
\end{tabbing}
\epf

Let us define
\begin{definition}
A filter $\F$ on $A$ is
{\em weakly flat} if and only if $$\Limsup_{\F} \Mm(\F) = 0.$$
\end{definition}
The previous definition may be interpreted as a generalization 
to non-symmetric spaces of the idea that 
the diameter of the elements of the filter may be chosen
arbitrary small.
Let us rephrase this definition. 
A filter $\F$ on $A$ is weakly flat if and only if 
for any $\epsilon > 0$, there exists an $f \in \F$
such that for any element $x$ of $f$, for any $g \in \F$, there exists 
$y \in g$ such that $A(x,y) \leq \epsilon$.\\

We shall introduce also the following filters
whose relevance will appear later.
\begin{definition}
A filter $\F$ on $A$ is {\em flat} if and only if 
for any $\epsilon > 0$, there exists an $f \in \F$
such that for any finite family of elements ${(x_i)}_{i \in I}$ of $f$, for any $g \in \F$, 
there exists $y \in g$ such that $A(x_i,y) \leq \epsilon$.\\
\end{definition}

A few remarks are in order.\\

One has the inclusion of classes of filters:
\begin{center} 
Cauchy $\Rightarrow$ flat $\Rightarrow$ weakly flat.
\end{center}
If the space $A$ is symmetric, that is when $A(x,y) = A(y,x)$,
then any flat filters on $A$ is also Cauchy. We shall see 
later \ref{symA} a few consequences of this fact.\\  

Also one might think to consider the filters $\F$ satisfying 
\begin{point}
\label{clumsy}
$$\Limsup_{\F} \Mp(\F) = 0.$$
\end{point}
These filters are actually useless for 
the study of non-symmetric spaces as for the obvious 
example non-symmetric 
space the base $\R$ itself, they do not correspond to 
``oriented neighborhoods''.
Consider any real $x$ and define its neighborhood filter
as generated by the family 
$\{ y \mid [y,x] \leq \epsilon \}$, 
where $\epsilon > 0$. This filter does not satisfy
\ref{clumsy}.
For the same reason it will occur that the class of  
filters defined by \ref{clumsy} cannot generally 
correspond to any class of modules containing 
the representable presheaves.
We have therefore focussed on the weakly flat filters 
and the associated operator $\Mm$.\\ 

The operator sending modules to the associated 
filters seems obvious. 
With $\wp(X)$ denoting the 
powerset of $X$ with inclusion ordering.
\begin{definition}
For any module $M$, Let $\Ba(M)$ denote the subset of $\wp(A)$, 
$\{ \Ba(M)(\epsilon) \mid \epsilon \in ]0, + \infty] \}$,
where $\Ba(M)(\epsilon)$ denotes the set $\{ x \in A \mid 
M(x) \leq \epsilon\}$. Let also $\F(M)$ denote the upper 
closure in $\wp\wp(A)$ of $\Ba(M)$. 
\end{definition}

We are ready to establish correspondences
between various ($\R$- !) categories of modules and filters
as well as a few other ``$\bullet$'' points.
Let $\FFu(A)$ and $\FFd(A)$ will stand for the
ordinary categories with objects respectively weakly flat 
filters, and flat filters on $A$, both with reverse inclusion ordering.
We are going to prove
\begin{theorem}
\label{flatref}
The map $\F \mapsto  \Mm(\F)$ determines
reflectors $\FFu(A) \rightarrow \FPu(A)_0$
and $\FFd(A) \rightarrow \FPd(A)_0$.
Their respective right adjoints
$\FPu(A)_0 \hookrightarrow \FFu(A)$
and $\FPd(A)_0 \hookrightarrow \FFd(A)$
send any module $M$ to the filter $\F(M)$ with
basis $\Ba(M)$. They are full.
\end{theorem}
Moreover the inclusions $\FPd(A)_0 \hookrightarrow \FPu(A)_0$
and $\FFd(A) \hookrightarrow \FFu(A)$ are maps between the above 
adjunctions.
This reflection will yield a
notion of morphisms between weakly flat filters more
general than the inclusion ordering. We shall later consider the 
associated category of fractions $\FFu^*(A)$ that is equivalent 
to $\FPu(A)_0$. In this category weakly flat filters and flat 
filters will be defined then in terms of colimits of the so-called 
{\em forward Cauchy sequences} (\ref{charffil}).\\    

{\it $\bullet$ $\FPu$-modules and weakly flat filters.}\\

We shall establish the reflection
$\Mm \dashv \F: \FFu(A) \hookleftarrow {\FPu(A)}_0$ 
as well as a couple of results regarding weakly flat 
filters. This full reflection results from 
\ref{OFl1}, \ref{OFl2}, \ref{counit}, \ref{unit}, 
\ref{wwf}, \ref{cot}, \ref{wfpt} and
\ref{counitiso} below.\\  

\begin{fact}
\label{OFl1}
For any modules $M_1$ and $M_2$ on $A$,
if $M_1 \Rightarrow M_2$ then
$\F(M_1) \supseteq \F(M_2)$.
\end{fact}

\begin{fact}
\label{OFl2}
For any filters $\F_1$ and $\F_2$ on $A$,
if $\F_1 \supseteq \F_2$  then
$\Mm(\F_1) \Rightarrow \Mm(\F_2)$.
\end{fact}

\begin{proposition}
\label{counit}
Let $M$ be a left module on $A$ then
for all $x$, $M(x) \leq \bigvee_{\epsilon > 0} \bigwedge_{y \mid M(y) \leq \epsilon} A(x,y)$,
i.e. 
$M \Leftarrow \Mm \circ \F(M)$.
\end{proposition}
\pf
Let $x \in A$.
For all $y \in A$, $M(x) \leq M(y) + A(x,y)$ thus
for all $y \in A$, such that $M(y) \leq \alpha$, $M(x) \leq A(x,y) + \alpha$
and $M(x) \leq \bigwedge_{y \mid M(y) \leq \alpha} A(x,y) + \alpha$.
Consider $\epsilon > 0$. 
The map $\alpha \mapsto \bigwedge_{y \mid M(y) \leq \alpha} A(x,y)$
reverses the order so $\bigwedge_{ y \mid M(y) \leq \epsilon } A(x,y)$
$=$ $\bigvee_{ \alpha \geq \epsilon} \bigwedge_{ y \mid M(y) \leq \alpha } A(x,y)$
and $(*)$
$M(x)$ $\leq$ $\bigvee_{ \alpha \geq \epsilon} \bigwedge_{ y \mid M(y) \leq \alpha } A(x,y) + \epsilon$.
Also for any $\alpha \leq \epsilon$,
\begin{tabbing}
\hspace{1cm}$M(x)$ \=$\leq$  \=$( \bigwedge_{ y \mid M(y) \leq \alpha } A(x,y) ) + \alpha$\\  
\>$\leq$  \>$( \bigwedge_{ y \mid M(y) \leq \alpha } A(x,y) ) + \epsilon$
\end{tabbing}
and thus $(**)$
$M(x)$ $\leq$  $\bigvee_{\alpha \leq \epsilon}  \bigwedge_{ y \mid M(y) \leq \alpha } A(x,y)  + \epsilon$.
$(*)$ and $(**)$ give
$M(x) \leq  \bigvee_{\alpha > 0} \bigwedge_{y \mid M(y) \leq \alpha}A(x,y) + \epsilon$.\\
\epf

\begin{proposition}
\label{unit}
A filter $\F$ on $A$ is weakly flat if and only if
$\F \supseteq  \F \circ \Mm ( \F )$. 
\end{proposition}
\pf One has the successive equivalences.
\begin{tabbing}
\hspace{0.5cm}\=$\F$ is weakly flat\\ 
\hspace{1cm}if and only if\\
\>$\bigwedge_{f \in \F} \bigvee_{x \in f} \Mm(\F) = 0$\\
\hspace{1cm}if and only if\\
\>for all $\epsilon > 0$, 
there exists $f \in \F$ such that 
for all $x \in f$, $\Mm(\F)(x) \leq \epsilon$,\\
\hspace{1cm}if and only if\\
\>for all $\epsilon > 0$, 
there exists $f \in \F$ such that 
$f \subseteq \Ba(\Mm(\F))(\epsilon)$\\
\hspace{1cm}if and only if\\
\>$\F \supseteq  \F \circ \Mm ( \F )$.
\end{tabbing}
\epf

Let us note at this stage
\begin{proposition}
\label{zoi}
For any weakly flat filter $\F$ on $A$ and any left module $M$ on $A$,
$\F \supseteq \F(M)$ if and only if
$\Mm(\F) \Rightarrow M$. 
\end{proposition}
\pf
If $\F \supseteq \F(M)$ then
$\Mm(\F) \Rightarrow \Mm \circ \F(M) \Rightarrow M$ by
\ref{OFl2} and \ref{counit}. Conversely
if $\Mm(\F) \Rightarrow M$ then
$\F \supseteq \F \circ \Mm(\F) \supseteq \F(M)$ 
by \ref{OFl1} and \ref{unit}.
\epf

\begin{proposition}
\label{fac22} For any right module $N : \xymatrix{ A \ar[r]|{\circ} &  I}$,
and any filter $\F$ in $A$,
$$N * \Mm(\F) \geq \Liminf_{\F} N$$
If $\F$ is moreover weakly flat then
the previous inequality becomes an equality.
\end{proposition}
\pf
For any module $N$ and any filter $\F$ as above,
\begin{tabbing}
\hspace{1cm}$N*\Mm(\F)$\=$=$\=$\bigwedge_{x \in A} ( \Mm(\F)(x) + N(x) )$\\
\>$=$ \>$\bigwedge_{x \in A} ( \bigvee_{f \in \F} \bigwedge_{y \in f} A(x,y)  )  + N(x)  )$\\
\>$\geq$ \>$\bigwedge_{x \in A} \bigvee_{f \in \F} ( (\bigwedge_{y \in f} A(x,y))    + N(x))$\\
\>$=$ \>$\bigwedge_{x \in A} \bigvee_{f \in \F} \bigwedge_{y \in f} ( A(x,y)    + N(x))$\\
\>$\geq$ \>$\bigvee_{f \in \F} ( \bigwedge_{y \in f} N(y) )$.
\end{tabbing}
Let us suppose moreover that $\F$ is weakly flat.
Let $\epsilon > 0$. One may choose $f_{\epsilon} \in \F$ such that
when $x \in f_{\epsilon}$, $\Mm(\F)(x) \leq \epsilon$.
Thus \begin{tabbing}
\hspace{1cm}$N * \Mm(\F)$\=$=$\=$\bigwedge_{x \in A} ( \Mm(\F)(x) + N(x) )$\\ 
\>$\leq$ \>$\Mm(\F)(x) + N(x)$, for any $x \in f_{\epsilon}$\\
\>$\leq$ \>$\epsilon + N(x)$, for any $x \in f_{\epsilon}$.
\end{tabbing}
Thus 
\begin{tabbing}
\hspace{1cm}$N * \Mm(\F)$\=$\leq$\=
$\bigwedge_{x \in f_{\epsilon}} ( \epsilon + N(x) )$\\
\>$=$ 
\>$\epsilon + \bigwedge_{x \in f_{\epsilon}} N(x)$\\
\>$\leq$ 
\>$\epsilon + \bigvee_{f \in \F} \bigwedge_{x \in f} N(x)$.     
\end{tabbing}
\epf

\begin{fact}
\label{term} 
For any module $M: \xymatrix{I \ar[r]|{\circ} & A}$ if 
$- * M: [A,\R] \rightarrow \R$ preserves 
the terminal object (i.e. $\bigwedge_{ x \in A} M(x) = 0$)
then $\F(M)$ is a filter on $A$ with basis the family $\Ba(M)$.
\end{fact}
\pf  
Let us see first that the set of subsets of the form
$\Ba(M)(\epsilon)$ for $\epsilon > 0$, is a filter
basis on $A$.
$\Ba(M)$ is trivially cofiltered subset of $\wp(A)$ ordered 
by inclusion.
Since $\bigwedge_{ x \in A} M(x) = 0$, for any $\epsilon > 0$
there is one $x$ with $M(x) \leq \epsilon$, i.e. 
$\Ba(M)(\epsilon) \neq \emptyset$.
\epf

As a consequence of \ref{term} and \ref{P1flatM} below,
one gets
\begin{corollary}
\label{wwf} 
For any module $M: \xymatrix{I \ar[r]|{\circ} & A}$ 
if $M$ is $\Pu$-flat then $\F(M)$ is weakly flat.  
\end{corollary}

\begin{lemma}
\label{P1flatM}
If $M: \xymatrix{I \ar[r]|{\circ} & A}$ is $\Pu$-flat
then for any $\epsilon > 0$ and any $x$ 
with $M(x) < \epsilon$ and any $\alpha > 0$,
there is a $y$ such $M(y) \leq \alpha$ and
$A(x,y) \leq \epsilon$.
\end{lemma}
\pf
$-*M$ preserves cotensors and the (conical) terminal object.
Consider $\epsilon > 0$ and $x$ with
$M(x) < \epsilon$. Then 
\begin{tabbing}
\hspace{1cm}$0$ \=$=$ \=$[M(x), M(x)]$\\
\>$=$ \>$[M(x), A(x,-) * M ]$\\
\>$=$ \>$( M(x) \pitchfork A(x,-) ) * M$\\
\>$=$ \>$\bigwedge_{y \in A} ( M(y) + [M(x), A(x,y)] )$.
\end{tabbing}
So for any $\delta > 0$, there is an $y$ such that 
$M(y) + [M(x), A(x,y)] \leq \delta$.
This $y$ satisfies $M(y) \leq \delta$,
and $A(x,y) \leq M(x) + \delta$.
Now given any $\alpha>0$,
considering $\delta = min \{ \alpha, \epsilon - M(x) \}$
one may find a $y$ as required.
\epf

$\R$ has a very peculiar property that we are going to use.
\begin{fact}
\label{facR}
For any $v$ in $\R$ and any non empty family $(a_i)_{i \in I}$ in $\R$,
$[v,\bigwedge_{i \in I} a_i] = \bigwedge_{i \in I} [v,a_i]$.
\end{fact}
\pf
Since $[v,-]$ preserves the usual ordering on $\R$,  
$[v, \bigwedge_{i \in I} a_i] \leq \bigwedge_{i \in I} [v, a_i]$
(even if $I$ is empty).
Conversely, fix $\epsilon > 0$.  
Since $I$ is not empty, there exists $j \in I$ such that
$\epsilon + \bigwedge_{i \in I} a_i \geq a_j$. Also
$[v, \bigwedge_{i \in I} a_i] \geq [v, \bigwedge_{i \in I} a_i]$,
so $v + [v, \bigwedge_{i \in I} a_i] \geq \bigwedge_{i \in I} a_i$.
For a $j$ as above,
$\epsilon + v + [v, \bigwedge_{i \in I} a_i] \geq a_j$
and $\epsilon + [v, \bigwedge_{i \in I} a_i] \geq [v, a_j] 
\geq \bigwedge_{i \in I}[v,a_i]$.  
\epf

\begin{fact}
\label{cot}
If $\F$ is a weakly flat filter on $A$ then 
$ - * \Mm(\F)$ preserves cotensors.
\end{fact}
\pf Given $v \in \R$ and $N: \xymatrix{ A \ar[r]|{\circ} & I}$
we have to show $(v \pitchfork N) * \Mm(\F) =  [v, N*\Mm(\F)]$. 
According to \ref{fac22}, 
$(v \pitchfork N) * \Mm(\F) = \bigvee_{f \in \F} \bigwedge_{x \in f} [v,N(x)]$
and $[v, N * \Mm(\F) ] = [v,  \bigvee_{f \in \F} \bigwedge_{x \in f} N(x)]$
$=$ $\bigvee_{f \in \F} [v, \bigwedge_{x \in f} N(x)]$.
Since all the $f \in \F$ are non empty the result follows then from \ref{facR}.
\epf

\begin{proposition}
\label{wfpt}
If $\F$ is a weakly flat filter on $A$
then $-* \Mm(\F)$ preserves the terminal object. 
\end{proposition}
\pf
We have to show that $\bigwedge_{x \in A} \Mm(\F)(x) = 0$.
For any $\epsilon > 0$, since $\Limsup_{\F}  \Mm(\F) = 0$
one may find an $f \in \F$ such that for any $x \in f$,
$\Mm(\F)(x) \leq \epsilon$. Since that $f$ is not empty
then $\bigwedge_{x \in A} \Mm(\F)(x) \leq \epsilon$
\epf

\begin{proposition}
\label{counitiso}
Let $M$ be a $\Pu$-flat left module on $A$ then
for all $x$, 
$M(x)$ $\geq$ $\bigvee_{\epsilon > 0} \bigwedge_{y \mid M(y) \leq \epsilon}A(x,y)$ $=$ 
$\Mm \circ \F(M)$.
\end{proposition}
\pf
Let $x \in A$ and $\epsilon$ be such that $M(x) <  \epsilon$. 
According to \ref{P1flatM},
for all $\alpha$, there exists $y$ such that $M(y) \leq \alpha$ and
$A(x,y) \leq \epsilon$. Thus for all $\alpha$, 
$\bigwedge_{y \mid M(y) \leq \alpha} A(x,y) \leq \epsilon$, i.e
$\bigvee_{\alpha > 0} \bigwedge_{y \mid M(y) \leq \alpha} A(x,y) \leq \epsilon$.
\epf

{\it $\bullet$ $\FPd$-modules and flat filters.}\\

Now we establish the reflection $\F: \FFu(A) \hookleftarrow {\FPu(A)}_0$.
This results from \ref{flatf} and \ref{lex} below.
\begin{fact}
\label{flatf} 
For any $\Pd$-flat module $M: \xymatrix{I \ar[r]|{\circ} & A}$,   
$\F(M)$ is a flat filter.
\end{fact}
\pf
If $- * M$ preserves conical finite limits then it preserves
in particular the terminal object 
and according to \ref{term},$\F(M)$ is a filter on $A$. 
The fact that the filter basis $\Ba(M)$ 
generates a flat filter is a consequence of the following lemma.
\epf

\begin{lemma}
\label{P2flatM}
If $M$ is $\Pd$-flat 
then for any $\epsilon > 0$ and
and any finite family ${(x_i)}_{i \in I}$ 
such that for all $i$, $M(x_i) < \epsilon$ and any $\alpha > 0$,
there is a $y$ such $M(y) \leq \alpha$ and for all $i \in I$,
$A(x_i,y) \leq \epsilon$.
\end{lemma}
\pf
$-*M$ preserves conical finite limits and cotensors.
Consider $\epsilon > 0$ and a finite family of $x_i$'s such that 
$M(x_i) < \epsilon$. 
Let us write $\epsilon' = \bigvee_{i \in I}M(x_i)$. Then 
\begin{tabbing}
\hspace{1cm}$0$ \=$=$ \=$[\epsilon', \bigvee_{i \in I}M(x_i)]$\\
\>$=$ 
\>$[\epsilon', \bigvee_{i \in I}( A(x_i,-) * M )]$\\
\>$=$
\>$[\epsilon', ( \bigvee_{i \in I}A(x_i,-) ) * M ]$\\ 
\>$=$$(  \epsilon' \pitchfork (\bigvee_{i \in I}A(x_i,-) ) * M$\\
\>$=$$\bigwedge_{y \in A} (M(y) + [\epsilon', \bigvee_{i \in I}A(x_i,y)])$.
\end{tabbing}
So for any $\delta > 0$, there is an $y$ such that 
$M(y) + [\epsilon', \bigvee_{i \in I}A(x_i,y)] \leq \delta$.
This $y$ satisfies $M(y) \leq \delta$,
and for all $i$, $A(x_i,y) \leq \epsilon' + \delta$.
Now given any $\alpha > 0$, by considering $\delta = min \{ \alpha, \epsilon - \epsilon' \}$
one may find a $y$ as required.
\epf

\begin{fact}
\label{lex}
If the filter $\F$ is flat then
$-* \Mm(\F)$ preserves conical finite limits, i.e.
for any finite family $(N_i)_{i \in I}$ of right modules
on $A$
$$\bigwedge_{x \in A} ( \Mm(\F)(x) + \bigvee_{i \in I} N_i(x)  ) =  
\bigvee_{i \in I} \bigwedge_{x \in A} ( \Mm(\F)(x) +  N_i(x) ).$$
\end{fact}
\pf We shall only prove 
$\bigwedge_{x \in A} ( \Mm(\F)(x) + \bigvee_{i \in I} N_i(x)  ) 
\leq  \bigvee_{i \in I} \bigwedge_{x \in A} ( \Mm(\F)(x) +  N_i(x) )$ since
the reverse inequality is trivial.\\

Let $\epsilon > 0$.\\
If there is a filter $\F$ on $A$ then $A$ is not empty and
for each $i \in I$, there is an $x_i \in A$ such that
$$N_i * \Mm(\F) + \epsilon 
= \bigwedge_{x \in A} (\Mm(\F)(x) + N_i(x)) + \epsilon \geq 
\Mm(\F)(x_i) + N_i(x_i).$$

Let $f \in \F$. Given a family
of $x_i$'s as above, for each $i$,
$\Mm(\F)(x_i) \geq \bigwedge_{y \in f} A(x_i,y)$, 
thus there is an $y_i \in f$ such that
$\Mm(\F)(x_i) + \epsilon \geq A(x_i, y_i)$ and
\begin{tabbing}
\hspace{1cm}
$2 \cdot \epsilon + N_i * \Mm(\F)$  \=$\geq$ 
\=$A(x_i, y_i) + N_i(x_i)$\\
\>$\geq$ \>$N_i(y_i)$.
\end{tabbing}

Since $\F$ is flat, we can choose $f$ so that  
for the $y_i \in f$ as above, for all $g \in \F$, there exists 
$z \in g$ such that for all $i$, $A(y_i,z) \leq \epsilon$.
Thus for all $g \in \F$, there exists $z \in g$ such that 
for all $i$, 
\begin{tabbing}
\hspace{1cm}
$3 \cdot \epsilon   + N_i * \Mm(\F)$
\=$\geq$   
\=$A(y_i,z) + N_i(y_i)$ for some suitable $y_i$'s,\\
\>$\geq$ \>$N_i(z)$.
\end{tabbing}

Because $\epsilon$ is arbitrary we have shown so far that for any $g \in \F$, 
$$\bigvee_{i \in I}( N_i * \Mm(\F) ) \geq \bigwedge_{z \in g} \bigvee_{i \in I} N_i(z).$$
So \begin{tabbing}
\hspace{1cm}$\bigvee_{i \in I} ( N_i * \Mm(\F) )$ \=$\geq$
\=$\bigvee_{g \in \F} \bigwedge_{z \in g} \bigvee_{i \in I} N_i(z)$\\
\>$=$ \>$( \bigvee_{i \in I} N_i ) * \Mm(\F)$, according to \ref{fac22}. 
\end{tabbing}
\epf

{\it $\bullet$ The right morphisms for weakly flat filters.}\\

For any weakly flat filter $\F$ on $A$, that 
$\F = \F \circ \Mm(\F)$ is to say that
for all $f \in \F$ there exists $\epsilon > 0$ such that
$\Ba ( \Mm ( \F ) ) (\epsilon) \subseteq f$.
Call a weakly flat filter {\em closed} when it satisfies the latter 
condition. 
$\CFFu(A)$ denoting the full subcategory of $\FFu(A)$
of closed weakly flat filters, 
the inclusion $\CFFu(A) \hookrightarrow \FFu(A)$
has left adjoint $\F \circ \Mm$ according to \ref{flatref}.
From this one gets a notion of morphisms of weakly flat filters
more relevant than the inclusion as it will yield the 
characterization  of weakly flat/flat filters in terms of forward 
Cauchy sequences \ref{charffil}.
\begin{definition}
Let $\F_1$ and $\F_2$ be weakly flat filters on $A$.
We write $\F_1 \rightarrow \F_2$ if and only if
for all $\epsilon > 0$, there exists $f \in \F_1$ such that
for all $x \in f$, for all $g \in \F_2$, there exists
$y \in g$ such that $A(x,y) \leq \epsilon$. 
\end{definition}
Note that $\F_1 \rightarrow \F_2$ is by definition
$\F_1 \supseteq \F \circ \Mm(\F_2)$.
Thus weakly flat filters on $A$ with the above relation $\rightarrow$ define
a preorder denoted $\FFu^*(A)$ equivalent to $\FPu(A)_0$
($\FFu^*(A)$ is the category of fractions induced
by the reflector $\FFu(A) \rightarrow \FPu(A)_0$ - see \cite{Bor94}
prop 5.3.1 p.190).
Equivalence classes in $\FFu^*(A)$ are in one
to one correspondence with closed weakly flat filters.
In the same way, $\CFFd(A)$ will denote the full sucategory
of $\FFu(A)$ with objects closed flat filters and
$\FFd^*(A)$ will denote the full subcategory of
$\FFu^*(A)$ with objects flat filters.
$\FFd^*(A)$ is equivalent to $\FPd(A)_0$.\\

Closed weakly flat filters play a similar role
for non-symmetric spaces as the minimal Cauchy filters 
do for symmetric spaces. Note that   
\begin{fact}
\label{symA}
If $A$ is symmetric, 
\begin{itemize}  
\item $(1)$ flat filters on $A$ are Cauchy;
\item $(2)$ For any Cauchy filter $\F$, $\Ml(\F) = \Mr(\F)$. 
\item $(3)$ Any left adjoint module on $A$ 
has the same underlying map as its right adjoint;
\item $(4)$ For any left adjoint module $M$ on $A$,
$\F(M) = \Fs(M)$;
\item $(5)$ $\Pd$-flat modules are left adjoint;
\item $(6)$ Closed Cauchy filters are exactly the minimal Cauchy filters;
\item $(7)$ $\CFFd(A)$ is discrete.
\end{itemize}
\end{fact}
\pf
$(1)$ is already known.\\
$(2)$ trivial.\\
$(3)$ For any left module $M$ with right adjoint $\tilde{M}$,
according to $(2)$ their 
underlying maps satisfy $M = \Ml \circ \Fs(M) = \Mr \circ \Fs(M) = \raM$.\\
$(4)$ straightforward from $(3)$.\\ 
$(5)$ One has the successive equivalences
\begin{tabbing}
\hspace{1cm}\=$M$ is $\Pd$-flat\\
\>if and only if
$M = \Mm(\F)$ for a flat filter (\ref{flatref})\\
\>if and only if
$M = \Mm(\F)$ for a Cauchy filter, ($1$)\\
\>if and only if
$M$ is left adjoint.
\end{tabbing}
$(6)$ A Cauchy filter $\F$ is closed if and only
$\F = \F \circ \Mm(\F) = \F^s \circ \Mm(\F)$    
if and only if it is minimal as a Cauchy filter.\\
$(7)$ Actually the underlying
subcategory ${\cal C}$ of the full subcategory of presheaves
$[A^{op},\R]$ with objects left adjoint modules is 
- in this particular case - discrete.
And $\CFFd(A)$ is equivalent to ${\cal C}$ according
to $(6)$. 
Let us see that ${\cal C}$ is discrete.
For any left adjoint module $M$ on $A$, 
$M$ has the same underlying map as
its right adjoint $\raM$.
Now consider another left adjoint module $N$ on $A$,
with right adjoint $\tilde{N}$.
Then $M \Rightarrow N$ if and only if
$\forall x \in A, M(x) \geq N(x)$ 
if and only if
$\forall x \in A, \tilde{M}(x) \geq \tilde{N}(x)$
if and only if
$\tilde{M} \Rightarrow \tilde{N}$.
But also if $M \Rightarrow N$ then 
$1 \Rightarrow \tilde{M}N$
since $M \dashv \tilde{M}$
and 
$\tilde{N} \Rightarrow \tilde{M}$
since $N \dashv \tilde{N}$.  
So $M \Rightarrow N$ if and only if $M=N$.
\epf

{\it $\bullet$ Forward Cauchy sequences.}\\

Given a sequence $(x_n)_{n \in \nit}$ on $A$, the associated filter,
still denoted $(x_n)$,   
has basis the family of sets $\{ x_p \mid p \geq n \}$.
We say that 
\begin{definition}
$(x_n)$ is:\\
- {\em weakly flat}, respectively {\em flat}, if the associated filter is so;\\
- {\em forward Cauchy} if and only
if $\forall \epsilon > 0, \exists N \in \nit, 
\forall m \geq n \geq N, A(x_n,x_m) \leq \epsilon$.
\end{definition} 
Any forward Cauchy sequence is obviously flat.
The relevance of forward Cauchy sequences for non-symmetric
spaces appears with following result:
\begin{theorem}
\label{charffil}
$\FFu^*(A)$ has the colimits of functors with non-empty domain.
The family of forward Cauchy sequences is dense
in $\FFu^*(A)$.
Flat filters are exactly the colimits
in $\FFu^*(A)$ of functors taking values
forward Cauchy sequences and with non-empty
filtered domain. 
\end{theorem}
This is proved below. We shall
explicit the colimits in $\FFu^*(A)$ of functors with non-empty domain 
in \ref{yogl3}. According to \ref{yogl3} and \ref{yogl1}, any weakly 
flat filter $\F$ 
is a canonical colimit in $\FFu^*(A)$ of a functor with non-empty
domain and taking values forward Cauchy sequences, 
moreover this colimit is filtered when $\F$ is flat according 
\ref{yogl2}.\\ 

To simplify our notation, for any 
$f \subseteq A$, any $\epsilon > 0$ and, any $F \subseteq \wp\wp(A)$, 
we let:\\
- $P(f,\epsilon, F)$ denote the property:
``for all $x$ in $f$, for 
all $g \in F$ there
exists $y \in g$ such that
$A(x,y) \leq \epsilon$'';\\
- $Q(f,\epsilon, F)$ denote the property:
``for any finite family $x_1, ...,x_n \in f$, for 
all $g \in F$ there
exists $y \in g$ such that for all $i \in \{1,...,n \}$,
$A(x_i,y) \leq \epsilon$''.\\
So to say that a filter $\F$ on $A$ is weakly flat (respectively.
flat) is to say that for all $\epsilon > 0$,
there exists $f \in \F$ such that 
$P(f,\epsilon, \F)$ (respectively.
$Q(f,\epsilon, \F)$). Also for weakly flat filters
$\F_1, \F_2$,
\begin{fact} $\F_1 \rightarrow \F_2$
if and only if for all $\epsilon > 0$, 
$\exists f \in \F_1, P(f,\epsilon, \F_2)$.
\end{fact}
When $\F_2$ is moreover flat, one has 
that 
\begin{fact}
\label{F2flat}
If $\F_1 \rightarrow \F_2$ then
for all $\epsilon > 0$, 
$\exists f \in \F_1, Q(f,\epsilon, \F_2)$.
\end{fact}
This is a consequence of the following lemma.
\begin{lemma}
\label{usef1}
Given $f \subseteq A$, $\epsilon > 0$ and
a flat filter $\F$, if 
$P(f , \epsilon, \F)$ then for all $\alpha > 0$,
$Q(f, \epsilon + \alpha , \F ).$ 
\end{lemma}
\pf
Let $\alpha > 0$. Since $\F$ is flat,
there is a $g_{\alpha} \in \F$ such that
$Q(g_{\alpha}, \alpha, \F)$. 
Consider a finite family $x_1,..., x_n$
in $f$.
Since $P(f , \epsilon, \F)$, there exist 
$y_1$,..., $y_n$ and in $g_{\alpha}$ such that 
$A(x_i,y_i) \leq \epsilon$ for all $i$. 
Since $Q(g_{\alpha}, \alpha, \F)$,
for any $g \in \F$, 
one may find $t$ in $g$ such that $A(y_i,t) \leq \alpha$
for all $i$, so that
$A(x_i,t) \leq \epsilon + \alpha$ for all $i$. 
\epf

Let $Fil(A)$ stand for the category of filters on $A$
with reverse inclusion ordering. $\FFu(A)$ and $\FFd(A)$ are
full subcategories of $Fil(A)$. Moreover,
\begin{proposition}
\label{yogl3} In $Fil(A)$,
any non empty family $I$ has a least upper bound, moreover
if $I$ is a family of weakly flat filters, its upper bound is weakly
flat.
\end{proposition}
\pf
Let $\F = \bigcap_{i \in I} \F_i$, where $I$ is 
a non-empty set. $\F$ is a filter.
Suppose first that all the $\F_i$'s are weakly flat.
Given $\epsilon > 0$, for any $i$ there exists $f_i \in \F_i$
such that $P(f_i,\epsilon, \F_i)$. For all these $f_i$,
$P(f_i, \epsilon, \F)$ since $\F \subseteq \F_i$. 
Thus $f = \bigcup_{i \in I} f_i$
that belongs to $\F$ satisfies $P(f, \epsilon, \F)$.\\  


Since $\CFFu(A)$ is full and reflective in $\FFu(A)$,
if $i$ denotes here the inclusion 
$\CFFu(A) \hookrightarrow \FFu(A)$, 
any functor $F:J \rightarrow \CFFu(A)$ 
will have a colimit if
the composite $i \circ F$ has a colimit.
So from \ref{yogl3}, $\CFFu(A)$ has colimits 
of functors with non-empty domains,
and since $\FFu^*(A) \cong \CFFu(A)$,  
$\FFu^*(A)$ has these colimits as well
(the colimit in $\CFFu(A)$ of any non-empty family $(\F_i)_{i \in I}$
is $\F \circ \Mm (\bigcap_{i \in I} \F_i)$).
Remark that one can prove directly that the ordinary category
${\FPu(A)}_0$ has colimits of functors with non-empty
domain: for any non-empty family $M_i$ of $\Pu$-flat
left modules on $A$, the pointwise $\bigwedge_{i \in I} M_i$ is 
$\Pu$-flat. Nevertheless this straightforward proof relies
on the non-categorical argument \ref{facR}.\\
{\it Note to the referee - that part may be omitted}\\
Given a non-empty family of $\Pu$-flat ${(M_i)}_{i \in I}$,
a right module $N$ and $v \in \R$
\begin{tabbing}
\hspace{1cm}$(\bigwedge_{i \in I} M_i) * (v \pitchfork N)$ \=$=$
\=$\bigwedge_{x \in A}(  (\bigwedge_{i \in I} M_i(x) ) + [v,N(x)] )$\\
\>$=$ \>$\bigwedge_{x \in A} \bigwedge_{i \in I} ( M_i(x) + [v,N(x)])$\\
\>$=$ \>$\bigwedge_{i \in I} \bigwedge_{x \in A} ( M_i(x) + [v,N(x)])$\\
\>$=$ \>$\bigwedge_{i \in I} ( M_i * (v \pitchfork N)  )$\\
\>$=$ \>$\bigwedge_{i \in I} [ v , (M_i* N)  ]$, 
since each $M_i * -$ preserves cotensors\\
\>$=$ \>$[ v , \bigwedge_{i \in I} ( M_i* N)  ]$, since $I$ is not empty \ref{facR},\\
\>$=$ \>$[v, \bigwedge_{i \in I}\bigwedge_{x \in A} ( M_i(x) + N(x) ) ]$,\\
\>$=$ \>$[v, \bigwedge_{x \in A}\bigwedge_{i \in I} ( M_i(x) + N(x) ) ]$,\\
\>$=$ \>$[v, \bigwedge_{x \in A}( ( \bigwedge_{i \in I} M_i(x) ) + N(x) ) ]$,\\
\>$=$ $[v, ( \bigwedge_{i \in I} M_i) * N ]$.
\end{tabbing}

\begin{proposition}
\label{yogl1}
Given a weakly flat filter $\F$ on $A$ and a left module $M$
such that $\Mm(\F) \not \Rightarrow M$
(i.e. $\Mm(\F) \not \geq M$), 
there is a forward Cauchy sequence 
$(y_n)$ such that $(y_n) \rightarrow \F$ and
$\Mm(y_n) \not \Rightarrow M$.
\end{proposition}
\pf 
By hypothesis 
there exists $x \in A$ such that 
$\bigvee_{f \in \F} \bigwedge_{y \in f} A(x,y) < M(x)$. 
Consider such an $x$. There exists $\alpha> 0$ 
such that for any $f \in \F$, there exists $y \in f$ such that 
$A(x,y) + \alpha < M(x)$. Note then that for such a $y$,
$A(x,y)$ is necessarely finite.\\

Since $\F$ is weakly flat, one can define
a sequence $(f_n)$ of elements of $\F$ such that for all $n \in \nit$, 
$f_{n+1} \subseteq f_n$ and $P(f_n , \alpha \cdot 2^{-2-n}, \F)$.
$(f_n)$ is defined inductively as follows.\\

Choose first $f_0$ such that $P(f_0, \alpha \cdot 2^{-2}, \F)$.\\

If $f_n$ is defined then 
one can find $g \in \F$ such that $P(g, \alpha \cdot 2^{-2-(n+1)}, \F)$
and let $f_{n+1} = f_n \cap g$.\\

Then one can build a sequence $(y_n)$ where
for all integer $n$, $y_n \in f_n$,
$y_0$ is such that $A(x,y_0) + \alpha < M(x)$,
and for all integer $n$, $y_n \in f_n$, $A(y_n,y_{n+1}) \leq \alpha \cdot 2^{-2-n}$.\\

Actually this ensures that:\\
$(1)$ $(y_n)$ is forward Cauchy;\\
$(2)$ $(y_n) \rightarrow \F$;\\
$(3)$ $\Mm(y_n) \not \Rightarrow M$.\\
  
$(1)$ holds since
for all $n \leq p \in \nit$,
\begin{tabbing}
\hspace{1cm}$A(y_n, y_p)$ \=$\leq$ \=$A(y_n, y_{n+1}) + ... + A(y_{p-1},y_p)$\\ 
\>$\leq$ \>$\alpha \cdot ( 2^{-2-n} + 2^{-2-(n+1)} + ... )$\\
\>$=$ \>$\alpha \cdot 2^{-1-n}$.
\end{tabbing}

$(2)$ holds since $(y_n)$ is forward Cauchy, 
for any $n \in \nit$, $\{ y_p / p \geq n \} \subseteq f_n$ and
$P(f_n, \alpha \cdot 2^{-2 -n},\F)$.\\ 

$(3)$ holds since
for all $n \in \nit$, 
\begin{tabbing}
\hspace{1cm}$A(x,y_n)$\=$\leq$\=$A(x,y_0) + A(y_0, y_1) + ... + A(y_{n-1}, y_{n})$\\
\>$\leq$\= $A(x,y_o) + \alpha/2$
\end{tabbing}
so $A(x,y_n) + \alpha/2 <  M(x)$.
Thus $\Mm(y_n)(x) < M(x)$.
\epf
Note that according to \ref{yogl1} any weakly flat filter $\F$ dominates
at least one forward Cauchy sequence as $lim_{\F} \Mm(\F) = 0$ and 
thus $\Mm(\F) \not \Rightarrow +\infty$ with $+\infty$ the constant
module with value $+\infty$.


\begin{proposition}
\label{yogl2}
If $(x_n)$ and $(y_n)$ are weakly flat sequences and 
$\F$ is a flat filter such that 
$(x_n) \rightarrow \F \leftarrow (y_n)$, then
there exists a forward Cauchy sequence $(z_n)$ such that
$(x_n) \rightarrow (z_n) \leftarrow (y_n)$
and $(z_n) \rightarrow \F$.  
\end{proposition}
\pf
Since $(x_n) \rightarrow \F$ and $\F$ is weakly flat, 
one can define a sequence of integers $(N_i)_{i \in \nit}$ such that 
$P(\{x_n \mid n \geq N_i\},2^{-i}, \F)$.\\

Define  for all integers $i$,
$X_i$ has the finite set 
$\{x_n / N_i \leq n < N_{i+1} \}$.
Note then that $P(X_i,2^{-i}, \F)$.\\

Analogously define a sequence  $(M_i)_{i \in \nit}$ such that
$P(\{y_n \mid n \geq M_i\},2^{-i}, \F)$, and the sets
$Y_i = \{y_n / M_i \leq n < M_{i+1} \}$.\\

One may also find for any integer $i$, 
a $f_i \in \F$ such that $P(f_i,2^{-i},\F)$.\\ 

We are going to build by recurrence a sequence $(z_n)$ such that, for all integer $i$:\\
$(1)$ for all $x \in X_i$, $A(x,z_{i+1}) \leq 2^{-i+1}$;\\ 
$(2)$ for all $y \in Y_i$, $A(y,z_{i+1}) \leq 2^{-i+1}$;\\
$(3)$ $z_i \in f_i$;\\  
$(4)$ $A(z_{i}, z_{i+1}) \leq 2^{-i+1}$.\\

Choose first $z_0 \in f_0$.\\   

Suppose that $z_{i} \in f_i$.
Since $P(X_i ,2^{-i}, \F)$, 
$P(Y_i ,2^{-i}, \F)$
and 
$P(f_i ,2^{-i}, \F)$ then
$P(X_i \cup Y_i \cup f_i ,2^{-i}, \F)$. And since $\F$ is flat,
according to \ref{usef1},
$Q(X_i \cup Y_i \cup f_i ,2^{-i+1}, \F)$.
Because $X_i$ and $Y_i$ are finite
one may find $z_{i+1} \in f_{i+1}$ satisfying the point 
$(1)$, $(2)$ and $(4)$ below.\\
 
According to $(4)$, $(z_n)$ is forward Cauchy. Also from $(1)$,
respectively $(2)$,
one deduces $(x_n) \rightarrow (z_n)$, respectively  $(y_n) \rightarrow (z_n)$.
According to $(3)$, for any $p \geq n \in \nit$, 
$\Mm(\F)(z_p) \leq 2^{-n}$ thus 
$(z_n) \supseteq \F \circ \Mm(\F)$. 
\epf

\subsection{Non symmetric completions of general metric spaces}
\label{mcomp}
Let $A$ denote a general metric space. We are going to 
explicit in topological/metric terms the $\Pu$-
and $\Pd$-completions of $A$.\\

The ordinary category $A_0$ is a preorder
with $x \rightarrow y$ if and only if
$0 \geq A(x,y)$.\\

$\FPu(A)$ and $\FPd(A)$ are the small $\R$-categories
with objects respectively the $\Pu$-flat presheaves
and $\Pd$-flat presheaves on $A$, and 
with homs given by  
$$Hom(F,G) = [A^{op},\R](F,G) = \bigvee_{a \in A} [ Fa, Ga ].$$
$A$ embeds fully and faithfully into $\FPu(A)$ (respectively 
into $\FPd(A)$) by 
$a \mapsto Ya = A(-,a)$. 
Alternatively, according to \ref{flatref}, one has the following 
metric/topological 
description. $\FPu(A)$, respectively $\FPd(A)$, is isomorphic 
to the general metric space with points closed weakly flat filters on $A$,
respectively closed flat filters, and 
with pseudo distance $d$ defined for any $\F_1$, $\F_2$ by
$$d(\F_1,\F_2) = [A^{op}, \R](\Mm(\F_1),\Mm(\F_2)).$$
Actually we shall give an expression of this distance that do not 
refer to presheaves.
First let us show  
\begin{proposition}
\label{dwflat2}
For any left module $M$ on $A$ and any filter $\F$,
$$[A^{op}, \R] ( \Mm(\F), M ) \leq \Limsup_{\F} M.$$
If $\F$ is weakly flat then the inequality above
becomes an equality.
\end{proposition}
\pf
To simplify notations, let $LHS$ and $RHS$ denote respectively  
$\bigvee_{x \in A} [\Mm(\F)(x),M(x)]$
and $\bigwedge_{f \in \F} \bigvee_{z \in f} M(z)$.\\

According to the definition of $\Mm(\F)$, for all $x \in A$,
for all $f \in \F$, $\Mm(\F)(x) \geq \bigwedge_{z \in f} A(x,z)$. 
So for any $x \in A$, $f \in \F$ and any $\epsilon >0$, there exists
a $z \in f$ such that $A(x,z) \leq \Mm(\F)(x) + \epsilon$.
For such a $z$, $M(x) \leq M(z) + A(x,z)$ and
$M(x) \leq M(z) + \Mm(\F)(x) + \epsilon$.
So 
\begin{tabbing}
\hspace{1cm}\=$\forall x \in A, \forall f \in \F, \forall \epsilon > 0, \exists z \in f$,
$[\Mm(\F)(x),M(x)] \leq M(z) + \epsilon$,\\
thus \>
$\forall x \in A, \forall f \in \F, \forall \epsilon > 0$,
$[\Mm(\F)(x),M(x)] \leq (\bigvee_{z \in f} M(z)) + \epsilon$,\\
thus \>
$\forall f \in \F, \forall \epsilon > 0$,
$LHS \leq 
(\bigvee_{z \in f} M(z)) + \epsilon$,\\
thus \>
$\forall f \in \F$, 
$LHS \leq \bigvee_{z \in f} M(z)$,\\
thus \>
$LHS \leq RHS$.\\
\end{tabbing}
Suppose now that $\F$ is weakly flat.
Consider $\epsilon > 0$. One may find an $f_{\epsilon} \in \F$ 
such that for all $z \in f_{\epsilon}$, 
$\Mm(\F)(z) \leq \epsilon$. 
So,
\begin{tabbing}
\hspace{1cm}\=for any $z \in f_{\epsilon}$, $M(z) \leq [\Mm(\F)(z), M(z)] +
\epsilon$,\\
thus \>
$\bigvee_{z \in f_{\epsilon}} M(z) \leq 
(\bigvee_{x \in A}  [\Mm(\F)(x), M(x)]) + \epsilon$,\\
and \>
$RHS \leq LHS + \epsilon$.
\end{tabbing}
\epf

According to the latter property, one gets
\begin{fact}
\label{wfHom}
For any weakly flat filters $\F_1$ and $\F_2$,
$$[A^{op}, \R](\Mm(\F_1), \Mm(\F_2)) =  
\Limsup_{x \in \F_1} \Liminf_{y \in \F_2} A(x,y).$$ 
\end{fact}

The limits in \ref{wfHom} ``commute'' when the first argument is Cauchy:
\begin{fact}
\label{CcHom}
For any Cauchy filter $\F_1$ and any weakly flat filter $\F_2$,
$$[A^{op}, \R](\Mm(\F_1),\Mm( \F_2)) = \Liminf_{y \in \F_2} \Limsup_{x \in \F_1} A(x,y).$$
\end{fact}
To see this we shall establish the following categorical result
\begin{proposition}
\label{Homladj}
Given a complete monoidal closed $\V$, 
if a $\V$-module  $M: \xymatrix{I \ar[r]|{\circ} & A}$
has a right adjoint $\tilde{M}: \xymatrix{A \ar[r]|{\circ} & I}$
then for any left $\V$-module $N$: 
$[A^{op}, \V](M,N) \cong \tilde{M} * N$.
\end{proposition}
\pf Given any presheaves $M,N: A^{op} \rightarrow \V$,
$[A^{op}, \V](M,N)$ is the right lifting of $N$ through 
$M$. Suppose that moreover $M$ has a right adjoint $\tilde{M}$. If
$\epsilon: M * \tilde{M} \rightarrow 1: \xymatrix{A \ar[r]|{\circ} & A}$ 
denotes the counit of the adjunction $M \dashv \tilde{M}$ then
it is easy to see that 
$M * ( \tilde{M}  * N) \cong \xymatrix{ ( M * \tilde{M} ) * N \ar[r]^{\epsilon * N} & N}$ 
exhibits $\tilde{M} * N$ as the right lifting of $N$ through $M$.
\epf

\pf (of \ref{CcHom})
Now remember that if $\F$ is a Cauchy filter on $A$ then $\Mm(\F) = \Mp(\F) = \Ml(\F)$ 
and this left module on $A$ has right adjoint the module 
$\Mr(\F)$ defined 
by the map $x \mapsto \Limsup_{y \in \F} A(y,x) = \Liminf_{y \in f} A(y,x)$.
So according to \ref{Homladj} and \ref{fac22}, for any Cauchy filter $\F_1$ and any 
weakly flat filter $\F_2$,
\begin{tabbing}
\hspace{1cm}$[A^{op}, \R]( \Mm(\F_1), \Mm(\F_2) )$ \=$=$  \=$\Mr(\F_1) * \Mm(\F_2)$\\
\>$=$
\>$\Liminf_{y \in \F_2} \Mr(\F_1)(y)$\\ 
\>$=$
\>$\Liminf_{y \in \F_2} \Limsup_{x \in \F_1} A(x,y).$
\end{tabbing}
\epf

One has a notion of non-symmetric convergence in $A$.
The {\em neighborhood filter} of $x \in A$, 
denoted $V_A(x)$, is 
the filter generated by the family of subsets  
$\{y \mid A(y,x) \leq \epsilon \}$ with $\epsilon > 0$.
Which is to say that $V_A(x)$ is $\F(A(-,x))$.
Given a filter $\F$ on $A$ and 
$x \in A$, we say that $\F$ {\em converges} to $x$,
that we write $\F \rightarrow x$, 
if and only if $\F \supseteq V_A(x)$.
If $\F$ is weakly flat then by \ref{zoi}  
$\F$ converges to $x$ if and only
if $\Mm(\F) \Rightarrow A(-,x)$.
By Yoneda, this is also equivalent to say 
that for any $a \in A$, 
$$A(x,a) \geq [A^{op},\R](\Mm(\F), A(-,a)),$$
or according to \ref{dwflat2} that 
$$A(x,a) \geq \Limsup_{\F} A(-,a).$$

It remains to explicit in topological terms
the notion of colimits indexed by flat presheaves.
\begin{definition}
A filter $\F$ on $A$ has {\em representative} $x_0$
if and only if for all $a \in A$,
$$A(x_0,a) = \Limsup_{\F}A(-,a).$$ 
\end{definition} 
Which is exactly to say that $x_0$ is the colimit $\Mm(\F)* 1$.
In particular if a representative of $\F$ exists
then it is unique up to isomorphism. In this case we denote it 
$rep(\F)$. Note that $rep(\F)$ when it exists is necessarely the 
greatest lower bound in $A_0$ amongst objects such that $\F$ 
converges to.\\

Given a filter $\F$ on $A$ and a map $G:A \rightarrow B$
the direct image of $\F$ denoted $G(\F)$ is the filter
on $B$ generated by the subsets $G(f)$ for $f \in \F$.
It is easy to check for $\F$ and $G$ as above that 
if $G$ is non-increasing
and $\F$ is weakly flat, respectively flat, 
then $G(\F)$ is again weakly flat, respectively flat. Moreover,
\begin{proposition}
\label{image}
Given a weakly flat filter $\F$ on $A$, and a functor $G: A \rightarrow B$,
$\Mm(G(\F)): B^{op} \rightarrow \R$ is the (pointwise) left Kan extension of 
$\Mm(\F): A^{op} \rightarrow \R$ along $G^{op}$.
\end{proposition}
\pf One has the pointwise computation (see \cite{Kel82},
(4.17), p.115), 
\begin{tabbing}
$Lan_{G^{op}}(\Mm(\F))(b)$ \=$=$ \=$B(b,G-) * \Mm(\F)$\\
\>$=$ \>$\bigvee_{f \in \F} \bigwedge_{x \in f} B(b,Gx)$, 
according to \ref{fac22},\\
\>$=$ \>$\bigvee_{g \in G(\F)} \bigwedge_{y \in g} B(b,y)$\\
\>$=$ \>$\Mm(G(\F))(b)$.
\end{tabbing}
\epf
Note that we could already infer from \ref{Lprpf}
that for any filter $\F$ and any functor $G: A \rightarrow B$,
if $\Mm(\F)$ is $\Pu$-flat, respectively 
$\Pd$-flat, then  $Lan_{G^{op}}(\Mm(\F))$
is also $\Pu$-flat, respectively $\Pd$-flat.\\

As a consequence of \ref{image},
\begin{corollary}
Given a weakly flat filter $\F$ on $A$, a functor $G: A \rightarrow B$ and
a presheaf $M : B^{op} \rightarrow \R$,
$$[B^{op},\R](\Mm(G(\F)),M) = [A^{op},\R](\Mm(\F), MG^{op})$$
\end{corollary}
And thus,
\begin{corollary} 
\label{trans}
Given a weakly flat filter $\F$ on $A$ and
a non expansive map $G: A \rightarrow B$,
for an object in $B$, 
to be the colimit $\Mm(\F) * G$ is 
equivalent to be representative of
the weakly flat filter $G(\F)$.
\end{corollary}
{\it note to the referee, to be omitted in the final version:}
\begin{tabbing}
\hspace{1cm}$B(\Mm(\F)*G,b)$\=$\cong$ \=$[A^{op},\R](\Mm(\F), B(G-,b))$\\
\>$\cong$ \>$[B^{op}, \R](\Mm(G(\F)), B(-,b))$.
\end{tabbing}
We shall call a general metric space $A$
{\em $\Pu$-complete}, respectively {\em $\Pd$-complete}
if the corresponding category $A$ is. 
So according to \ref{trans}, $A$ is $\Pu$-complete, respectively 
$\Pd$-complete if and only if
any weakly flat, respectively flat, filter on $A$ admits 
a representative.\\

Now given a weakly flat filter $\F$ on $K$, 
and non-increasing maps $G:K \rightarrow A$, 
and $H:A \rightarrow B$, $H$ (as a functor) preserves the 
colimit $\Mm(\F)*G$ if and only $H$ (as a non-expansive map) 
preserves the representative of $G(\F)$, i.e.
$$H( rep(G(\F)) ) = rep (H \circ G( \F  ) ).$$
To sum up: the $\R$-functors preserving $\FPu$-colimits 
(respectively $\FPd$-colimits) are exactly 
the non-expansive maps preserving the representatives of
weakly flat filters (respectively those of flat 
filters).\\

A direct translation of \ref{freecoc} gives
for any general metric space $A$ two completions. 
\begin{theorem}
For any general metric space $A$, there exists 
a $\Pu$-complete (respectively $\Pd$-complete) metric space 
$\bar{A}$ together with a map $i_A: A \rightarrow \bar{A}$ such that 
for any non-expansive $f: A \rightarrow B$ with codomain 
$B$ $\Pu$-complete (respectively $\Pd$-complete) there exists 
a unique $\bar{f}: \bar{A} \rightarrow B$ preserving representatives
(respectively representatives of flat filters) and such that 
$\bar{f} \circ i_a = f$. 
\end{theorem} 
For an $f: A \rightarrow B$ as above, 
if one considers the completion $\bar{A}$ as a space of filters on $A$, 
then the extension $\bar{f}$ sends any filter $\F$ to
the representative of its direct image by $f$ in $B$.
To check this just come back to the categorical formulation.
From \cite{Kel82} Theorem 4.97,  
$\bar{f}$ is the left Kan extension of $f$ along $i_A$
and sends any $M$ in $\Fp(A)$ to $M*f$, ($\p = \Pu,\Pd$).
Translate then with \ref{trans}.\\

The rest of this section investigates examples of these
two completions.\\ 

Recall from \cite{Kel82} (3.74)
\begin{proposition}
Any monoidal closed $\V$ that is complete as an ordinary category
is complete and cocomplete as a $\V$-category, i.e. $\V$ admits all 
limits and colimits indexed by small $\V$-categories.
\end{proposition}
So
\begin{corollary}
$\R$ is $\Pu$-complete and $\Pd$-complete.
\end{corollary}

We shall also show that
\begin{proposition}
The $\Pu$- and $\Pd$-completion of $\R$
are both isomorphic to $\R$.
\end{proposition}
This results from \ref{wfeqf},
\ref{idflat1}, \ref{idflat2} and \ref{HomwfR}
below.

\begin{fact}
\label{wfeqf}
Weakly flat filters on $\R$ are flat.
\end{fact} 
This results from the following lemma.
\begin{lemma} Let $\F$ be a weakly flat filter on $\R$.
Let $\epsilon > 0$, and $f \in \F$ such that
$P(f,\epsilon,\F)$ then 
$Q(f,\epsilon,\F)$.
\end{lemma}
\pf
Consider any finite family $x_1$, ..., $x_n$ in $f$. 
Let $g \in \F$. There exist $y_1$,...,$y_n$
in $g$ such that $[x_i, y_i] \leq \epsilon$, i.e. 
$y_i \leq x_i + \epsilon$,  
for $i = 1,...,n$. Choosing the least of those $y_i$'s, say 
$z$ one has $[x_i,z] \leq \epsilon$ for all $i$'s. 
\epf

Let us identify now the weakly flat filters on $\R$.
For a filter $\F$ on $\R$,
we write $\liminf(\F)$ for $\Liminf_{\F} id$ $=$ 
$\bigvee_{f \in \F} \bigwedge_{x \in f} x$.
\begin{fact}
\label{idflat1}
If $\liminf(\F) \neq \infty$ then $\F$ is
weakly flat.
\end{fact} 
\pf
Let $\epsilon > 0$. Since $\liminf(\F)$ is finite, one may 
consider an $f \in \F$ such that 
$\bigwedge f$ ($= \bigwedge_{x \in f}x$) $\geq$ $\liminf(\F) - \epsilon$ 
and thus for any $x \in f$, $\liminf(\F) \leq x + \epsilon$. 
Now given any $g \in \F$, there exists $y \in g$
such that $y \leq \liminf(\F) + \epsilon$.
For this $y$, for any $x \in f$, $y \leq x + 2 \cdot \epsilon$, i.e 
$[x,y] \leq 2 \cdot  \epsilon$.
\epf

\begin{fact}
\label{idflat2}
If $\liminf(\F) = \infty$ then there are two cases
either $\F$ is the principal filter generated by
$\infty$ or not.
In the first case $\F$ is weakly flat, 
in the second case $\F$ is not weakly flat.
\end{fact}
The first case is trivial: $\F$ is a neighborhood filter.
For the second case, consider $\epsilon > 0$ and $f \in \F$. 
Then $\bigwedge f  \neq \infty$ and
there exists $g \in \F$ such that
$\bigwedge f + 2 \cdot \epsilon < \bigwedge g$.
So one may find $x \in f$ such that for any 
$y \in g$, $x + \epsilon < y$, i.e. $[x,y] > \epsilon$.
\epf

Eventually, for weakly flat filters $\F_1$ and $\F_2$ on $\R$,
one has the successive equations:
\begin{fact}
\label{HomwfR}
\begin{tabbing}
\hspace{1cm}
$[\R,\R](\Mm(\F_1), \Mm(\F_2))$ 
\=$=$ \=$\Limsup_{x \in \F_1} \Liminf_{y \in \F_2} [x,y]$,\\
\>$=$ \>$\Limsup_{x \in \F_1} (\bigvee_{g \in \F_2} \bigwedge_{y \in g}[x,y])$,\\
\>$=$ \>$\Limsup_{x \in \F_1} (\bigvee_{g \in \F_2} [x,\bigwedge g])$,\hspace{1cm}\=$(1)$\\
\>$=$ \>$\Limsup_{x \in \F_1} [x, \bigvee_{g \in \F_2} \bigwedge g])$, \>$(2)$\\
\>$=$ \>$\Limsup_{x \in \F_1}[x, \liminf(\F_2)]$,\\
\>$=$ \>$\bigwedge_{f \in \F_1} \bigvee_{x \in f}[x, \liminf(\F_2)]$,\\
\>$=$ \>$\bigwedge_{f \in \F_1} [\bigwedge f, \liminf(\F_2)]$, \>$(3)$\\
\>$=$ \>$[\bigvee_{f \in \F_1} \bigwedge f, \liminf(\F_2)]$, \>$(4)$\\
\>$=$ \>$[\liminf(\F_1), \liminf(\F_2)]$.
\end{tabbing}
\end{fact}
\begin{itemize}
\item
$(1)$ holds due to \ref{facR} since
any $g \in \F_2$ is non empty;
\item $(2)$ holds since representables
$[x,-]$ preserves limits;
\item $(3)$ holds since representables
$[-,y]$ turns colimits into limits; 
\item actually $(4)$ holds due to the characterization of weakly 
flat filters on $\R$ \ref{idflat2} and
another peculiar property of $\R$, \ref{facR2} above.
\end{itemize} 
\begin{fact}
\label{facR2}
Given any $v$ in $\R$, 
and any non-empty family $(a_i)_{i \in I}$
in $\R$ that satisfies the condition that
if $\bigvee_{i \in I} a_i = + \infty$ then 
there exists at least one $j \in I$ such that
$a_j = +\infty$, then 
$$[\bigvee_{i \in I} a_i,v] = \bigwedge_{i \in I}[a_i,v].$$
\end{fact}
\pf
Since $[-,v]$ reverses the usual ordering on $\R$, 
certainly 
$[\bigvee_{i \in I} a_i, v] \leq \bigwedge_{i \in I} [a_i,v]$
(even if $I$ is empty).
Conversely, let us fix $\epsilon > 0$.
By assumption one may find $j \in I$ such that
$\bigvee_{i \in I} a_i \leq a_j + \epsilon$.
For such a $j$,
$[\bigvee_{i \in I}a_i,v] \geq [a_j + \epsilon,v] = [\epsilon,[a_j,v]]$,
so 
 $[\bigvee_{i \in I}a_i,v] + \epsilon \geq [a_j,v] \geq 
\bigwedge_{i \in I}[a_i,v]$.
\epf


Eventually we shall study the completions of symmetric 
spaces. For a symmetric general metric space $A$, 
its $\Pd$-completion is its Cauchy completion (\ref{symA}). 
But the $\Pu$-completion of $A$ may not be symmetric as shown below.
\begin{proposition}
\label{sycomp}
The $\Pu$-completion of a symmetric $A$ is the set 
of non-empty closed subsets in its Cauchy-completion $\bar{A}$
with pseudo distance $d$ given by 
$d(X,Y) = \bigvee_{x \in X} \bigwedge_{y \in Y} \bar{A}(x,y)$. 
\end{proposition}
To prove this we shall establish first
\begin{lemma}
\label{liminfF}
For any filter $\F$ on $A$, any set 
$X$ of filters such that $\F$ is 
the intersection of the filters in $X$ - that we write $\F =
\bigcap X$ - and any map $t: Obj(A) \rightarrow \R$, 
$$\Liminf_{\F} t 
= \bigwedge_{\varphi \in X} (\Liminf_{\varphi} t).$$
\end{lemma}
\pf
Let us note $m = \bigwedge_{\varphi \in X} (\Liminf_{\varphi} t)$.
For any $\varphi \in X$, $\varphi \supseteq \F$ so
$\Liminf_{\varphi} t \geq \Liminf_{\F} t$ and
$m \geq \Liminf_{\F} t$.\\
Conversely, let us consider any positive real $v \leq m$.
Then for any $\varphi \in X$, 
$v \leq \bigvee_{g \in \varphi} \bigwedge_{x \in g} t(x)$.
Fix $\epsilon > 0$. For any $\varphi \in X$, there exists
$g_{\varphi} \in \varphi$ such that
$v \leq \bigwedge_{x \in g_{\varphi}} t(x) + \epsilon$.
Let $f = \bigcup_{\varphi \in X} g_{\varphi}$. Then $f \in \F$, 
$\bigwedge_{x \in f} t(x) = 
\bigwedge_{\varphi \in X} \bigwedge_{x \in g_{\varphi}} t(x)$
and $v \leq \bigwedge_{x \in f} t(x) + \epsilon
\leq \Liminf_{\F}t + \epsilon$.   
\epf

\pf (of \ref{sycomp})
Let $\bar{A}$ denote the Cauchy completion 
of $A$, it is isomorphic to the 
metric space with objects closed Cauchy filters on $A$ with
distance given for all $\varphi$, $\psi$ by
$\bar{A}(\varphi, \psi)$ $=$ $\Mr(\varphi) * \Mm(\psi)$
$=$ $\Liminf_{\psi} \Mr(\varphi)$ 
by \ref{Homladj}/\ref{CcHom} and \ref{fac22}. 
Since $A$ is symmetric, $\bar{A}$ is symmetric, also
forward Cauchy sequences in $A$ are Cauchy.\\

Consider a closed weakly flat filter $\F$ on $A$. $\bar{\F}$ 
will denote the set of closed Cauchy filters containing $\F$.
According to \ref{symA}-(6) and \ref{charffil}, $\bar{\F}$ is 
not empty and $\F$ is the intersection of the filters in $\bar{\F}$.\\

We shall show now the following property:\\
$(*)$ For any subset 
$X$ of $\bar{A}$ such that $\F$ is the intersection 
of the filters in $X$ and any closed Cauchy 
filter $\varphi$, $$[A^{op}, \R](\Mm(\varphi),\Mm(\F)) 
= \bigwedge_{\psi \in X} \bar{A}(\varphi,\psi).$$

Consider a Cauchy filter $\varphi$.
Then
\begin{tabbing} 
\hspace{1cm}$[A^{op}, \R]( \Mm(\varphi), \Mm(\F) )$ \=$=$ 
\=$\Mr(\varphi) * \Mm(\F)$, since $\Mm(\varphi)$ is left
adjoint \ref{Homladj}/\ref{CcHom},\\
\>$=$
\>$\Liminf_{\F} \Mr(\varphi)$,
according to \ref{fac22},\\
\>$=$
\>$\bigwedge_{\psi \in X} \Liminf_{\psi} \Mr(\varphi)$, 
according
to \ref{liminfF},\\
\>$=$
\>$\bigwedge_{\psi \in X} \Mr(\varphi) * \Mm(\psi)$\\
\>$=$
\>$\bigwedge_{\psi \in X} \bar{A}(\varphi, \psi)$.
\end{tabbing}

As a consequence of $(*)$, for any subset $X$ of $\bar{A}$
such that $\F = \bigcap X$, the adherence $\bar{X}$ of $X$ 
in $\bar{A}$ is $\bar{\F}$.
Hence $\bar{\F}$ is the only closed subset $X$ in the metric 
space $\bar{A}$, such that $\F = \bigcap X$.\\   

Now given two closed weakly flat filters 
$\F_1$ and $\F_2$ on $A$,
\begin{tabbing}
$[A^{op}, \R](\Mm(\F_1),\Mm(\F_2))$
\=$=$ \=$[A^{op}, \R](\bigwedge_{\varphi \in \bar{\F_1}} \Mm(\varphi),\Mm(\F_2))$, according to \ref{liminfF},\\
\>$=$ \>$\bigvee_{\varphi \in \bar{\F_1}}[A^{op}, \R](\Mm(\varphi),
\Mm(\F_2))$,\\
\>$=$ \>$\bigvee_{ \varphi \in \bar{\F_1} } \bigwedge_{\psi \in \bar{\F_2}}
\bar{A}(\varphi,\psi)$, 
according to $(*)$ above.
\end{tabbing}
\epf

\begin{remark}
Consider an embedding $i_A: A \hookrightarrow \bar{A}$ of 
a general metric space $A$ into its Cauchy-completion.
For any subset $X$ of $A$, the closure in $\bar{A}$ of the
direct image $i_A(X)$ is the set of Cauchy filters adherent to 
$X$. Considering now the inverse image by $i_A$ of this set   
one obtains the closure of $X$ in $A$. Nevertheless the direct image 
by $i_A$ does NOT yield in general a surjective correspondence
from closed subsets of $A$ to closed subsets of $\bar{A}$.
A counter-example for this would be the real line with the 
usual metric less one point. The Cauchy completion just 
adds the missing point say $x$, $\{x\}$ is closed in the completion
but its inverse image by the canonical embedding is empty.
\end{remark}

\end{section}

\begin{section}{The case $\V = \B$.}
Preorders as enrichments over the ``boolean'' category $\B$
and their Cauchy-completion were treated in \cite{Law73}.
We shall recall briefly these results. 
Then we characterize in this context 
the $\Pu$- and $\Pd$- flat presheaves and
the $\Pu$-completion and show that the $\Pd$-completion 
coincide with the classic ``algebraic'' (or ``ideal'') completion.\\

\label{preos}
$\B$ stands for the two-object category generated by the graph
$\xymatrix{ 0 \ar[r] & 1}$. It is a partial order and it has 
a monoidal structure with tensor $\wedge$
(the logical ``and'') and unit 
$1$. $\B$ is closed as for all $x,y,z \in \B$,
$$x \wedge y \leq z \Leftrightarrow x \leq (y \Rightarrow z)$$  
where $\Rightarrow$ denotes the usual entailment relation.\\

$\B$-categories are just preorders: for any $\B$-category $A$, 
its associated preorder is defined by $x \rightarrow y$ if and only 
if $A(x,y) = 1$. Along the same line there are one-to-one correspondences
between
\begin{itemize}
\item $\B$-functors and order preserving maps;
\item Right modules on a $\B$-category $A$ 
and downward closed subsets, or ``downsets'', on the preorder $A$;
\item left modules and upper closed subsets, in a dual way.
\end{itemize}
Under the above correspondences the Lawvere Cauchy-completion
for $\B$-categories is the Dedekind-Mac Neille completion 
for preorders. Also the $\emptyset$-completion occurs as the 
so-called the {\em downward completion} that is defined for 
a preorder as the set of its downsets with inclusion ordering.\\

Let us focus on the $\Pu$- and $\Pd$-flatness.
Further on 
$A$ denotes a $\B$-category that we might freely see as a 
preorder.
Using \ref{cotpoint}, \ref{precon} and \ref{colimpre} again,
one gets
\begin{fact}
For any module $M: \xymatrix{I \ar[r]|{\circ} & A}$,
\begin{itemize}
\item $-*M: [A,\B] \rightarrow \B$ preserves
the (conical) terminal object i.e.
$1 * M = 1$, if and only if 
$$(1)\;\;\bigvee_{x \in A} M(x) = 1;$$  
\item $-*M$ preserves conical finite limits if and only\\ 
$(2)$ For any finite family of right modules 
$N_i: \xymatrix{A \ar[r]|{\circ} & I}$, $i$ ranging in $I$,
$$\bigvee_{x \in A}( M(x) \wedge \bigwedge_{i \in I} N_i(x) )
= \bigwedge_{i \in I}( \bigvee_{x \in A} M(x) \wedge N_i(x)  );$$
\item $-*M$ preserves cotensors if and only if\\
$(3)$ For any $v \in \B$ and 
any right module $N: \xymatrix{A \ar[r]|{\circ} & I}$, 
$$\bigvee_{x \in A}( M(x) \wedge  (v \Rightarrow N(x)) ) 
\;\;=\;\; (\; v \Rightarrow  \bigvee_{x \in A}( M(x) \wedge N(x) )\; ).$$
\end{itemize}
\end{fact}
Condition $(1)$ above is equivalent to the fact that 
$\I_M$, the downset corresponding to $M$,
is not empty.
Condition $(3)$ reduces for $v=1$ to the trivial equation $N*M = N*M$. 
For $v= 0$, it reduces to $1 = \bigvee_{x \in A} M(x)$,
that is $(1)$ again.
Eventually condition $(2)$ is equivalent to the fact that 
$\I_M$ is directed i.e.
any finite family in $\I_M$ has an upper bound in $\I_M$.
So there are bijections between
the following objects on $A$: 
\begin{itemize}
\item $\Pu$-flat left modules and non-empty downsets, 
\item $\Pd$-flat left modules and non-empty directed downsets.
\end{itemize}
From these observations, one obtains straightforwardly that
\begin{itemize}
\item $\FPu(A)$ as a preorder is the set of  
non-empty downsets on $A$ with inclusion ordering; 
\item $\FPd(A)$ as a preorder is $Alg(A)$, the {\em algebraic completion}
(or {\em ideal completion}) of $A$, that is the set of its 
non-empty directed downsets
with inclusion ordering.
\end{itemize}

Eventually given $M: A^{op} \rightarrow \B$ and $G: A \rightarrow B$,
that $b \in B$ is the colimit $M*G$ is equivalent to the fact 
that $b$ is the least upper bound in the preorder $B$ of the downset 
generated by the direct image of $\I_M$ by $G$.
So from \ref{pflatchar} one gets two completions:
\begin{theorem}
Given a preorder $A$, $\FPu(A)$ together with the order 
preserving map $i_A: A \rightarrow \FPu(A)$ sending $a$ to the 
downset generated by $a$, are such that for any
preorder $B$ with least upper bounds, for any 
non-empty set and any order preserving
map $f: A \rightarrow B$, there is a unique
$\bar{f}: \FPu(A) \rightarrow B$ preserving least upper 
bounds of non empty sets and satisfying
$\bar{f} \circ i_A = f$.
\end{theorem}
And also
\begin{theorem}  
Given a preorder $A$, $Alg(A)$ together with the order preserving 
map $j_A: A \rightarrow Alg(A)$ sending $a$ to the downset 
generated by $a$, are such that
for any preorder $B$ with least upper bounds
for any non-empty directed set and any order 
preserving map $f: A \rightarrow B$, there is a unique 
$\bar{f}: Alg(A) \rightarrow B$ preserving least upper bounds
of non-empty directed sets and satisfying
$\bar{f} \circ j_A = f$.
\end{theorem}
It is common in the literature to define 
directed subsets in a partial order
as non-empty. Partial orders
with all least upper bounds for non-empty 
directed sets are usually called {\em directed
complete} partial orders or {\em dcpo}'s.
Also monotone maps between dcpos that preserve
least upper bounds of non-empty directed 
subsets are called {\em continuous}.    
\end{section}

\bibliographystyle{alpha}

\end{document}